\documentclass[1p]{elsarticle}

\newdefinition{rmk}{Remark}

\usepackage{color}

\usepackage{graphicx,amssymb,lineno}
\usepackage{amsmath}

\usepackage{hyperref}

\usepackage{appendix}

\usepackage[latin1]{inputenc}
\begin{document}

\begin{frontmatter}



\title{Asymptotic analysis of a thin fluid layer flow between two moving surfaces \tnoteref{min}}
\tnotetext[min]{This work has been partially supported by Ministerio de Econom\'{\i}a y Competitividad of Spain, under grant MTM2016-78718-P with the participation of FEDER, and the European Union's Horizon 2020 Research and Innovation Programme,
under the Marie Sklodowska-Curie Grant Agreement No 823731 CONMECH.}


\author[JM]{J. M. Rodr\'{\i}guez},  \ead{jose.rodriguez.seijo@udc.es}
\address[JM]{Department of Mathematics, Higher Technical University College of Architecture, Universidade da Coru\~na, Campus da Zapateira, 15071 - A Coru\~na, Spain}

\author[R]{R.
Taboada-V\'azquez\corref{cor}}\ead{raquel.taboada@udc.es}
\address[R]{Department of Mathematics, School of Civil Engineering, Universidade da Coru\~na, Campus de Elvi\~na, 15071 - A Coru\~na,  Spain}

\cortext[cor]{Corresponding author}

\begin{abstract}

In this paper we study the behavior of an incompressible viscous fluid moving between two very close surfaces also in motion. Using the asymptotic expansion method we formally justify two models, a lubrication model and a shallow water model, depending on the boundary conditions imposed. Finally, we discuss under what conditions each of the models would be applicable.

\end{abstract}

\end{frontmatter}


\section{Introduction}
\label{sec-introduccion}

The asymptotic analysis method is a mathematical tool that has been widely used to obtain and justify reduced models, both in solid and fluid mechanics, when one or two of the dimensions of the domain in which the model is formulated are much smaller than the others. 

After the pioneering works of Friedrichs, Dressler and Goldenveizer (see \cite{FD1961} and \cite{Golden1962}) the asymptotic development technique has been used successfully to justify beam, plate and shell theories (see, for example, \cite{Rigolot1972}, \cite{CiarletDestuynder1979a}, \cite{CiarletDestuynder1979b}, \cite{Ciarlet1980}, \cite{BermudezViano}, \cite{AgTutekN}, and many others).

This same technique has also been used in fluid mechanics to justify various types of models, such as lubrication models, shallow water models, tube flow models, etc. (see, for example, \cite{Elrod}, \cite{Dridi}, \cite{Cimatti1983}, \cite{BayadaChambat1986}, \cite{Nazarov}, \cite{Zeytounian}, 
\cite{MTZ97}, \cite{Grenier}, \cite{BessonLaydi}, 
\cite{AzeradGuillen}, \cite{GerbeauPerthame}, \cite{FerrariSaleri}, \cite{Hu}, \cite{Marche}, \cite{BreschNoble}, \cite{DBMS2009}, 
\cite{RodTab1}-\cite{RodTab6}, 
\cite{Dean1927}, \cite{Dean1928}, 
\cite{MaPa2001}, \cite{MaPaPa2007}, 
\cite{Panasenko2012}, \cite{Panasenko2015a}, \cite{Panasenko2015b},
\cite{Gonzalo2016}, \cite{Gonzalo2019}, and many others).

In this work, we are interested in justifying, again using the asymptotic development technique, a lubrication model in a thin domain with curved mean surface. Following the steps of \cite{BayadaChambat1986}, but with a different starting point, we devote sections \ref{sec-domain} and \ref{sec-lubrication} to this justification. During the above process we have observed that, depending on the boundary conditions, other models can be obtained, which we show in section \ref{sec-thin-layer}. In this section we derive a shallow water model changing the boundary conditions that we had imposed in section \ref{sec-lubrication}: instead of assuming that we know the velocities on the upper and lower boundaries of the domain, we assume that we know the tractions on these upper and lower boundaries.

Thus, two new models are presented in sections 3 and 4 of this article. These models can not be found in the literature,  as far as we know. In addition, the method used to justify them allows us to answer the question of when each of them is applicable. In section \ref{conc} we discuss the models yielded, as well as the difference between one model and another depending on the boundary conditions, reaching the conclusion that the magnitude of the pressure differences at the lateral boundary of the domain is key when deciding which of the two models best describes the fluid behavior.

\section{Derivation of the model} \label{sec-domain}
\subsection{Original domain}
Let us consider a three-dimensional thin domain, $\Omega^{\varepsilon}_t$, 
filled by a fluid, that varies with time $t \in [0, T]$, given by
\begin{eqnarray}\Omega^{\varepsilon}_t&=&\left\{ (x_1^{\varepsilon},x_2^{\varepsilon},x_3^{\varepsilon})\in{R}^3:x_i(\xi_1,\xi_2,t)\leq x_i^{\varepsilon} \leq
x_i(\xi_1,\xi_2,t)+ h^\varepsilon(\xi_1,\xi_2,t)N_i(\xi_1,\xi_2,t),  \nonumber\right.\\
&&\left. (i=1,2,3), \ (\xi_1,\xi_2)\in D\subset \mathbb{R}^2
\right\} \label{eq-o-domain} \end{eqnarray} where
$\vec{X}_t(\xi_1,\xi_2)=\vec{X}(\xi_1,\xi_2,t)=(x_1(\xi_1,\xi_2,t),x_2(\xi_1,\xi_2,t),
x_3(\xi_1,\xi_2,t))$ is the lower
bound surface parametrization, $h^\varepsilon(\xi_1,\xi_2,t)$ is the gap between the two surfaces in
motion, and $\vec{N}(\xi_1,\xi_2,t)$ is the unit normal vector:
\begin{equation}
\vec{N}(\xi_1,\xi_2,t)=\dfrac{ \dfrac{\partial\vec{X}}{\partial
\xi_1} \times \dfrac{\partial\vec{X}}{\partial \xi_2}
}{\left\|\dfrac{\partial\vec{X}}{\partial \xi_1} \times
\dfrac{\partial\vec{X}}{\partial \xi_2}\right\|}
\end{equation}

The lower bound surface is assumed to be regular, 
\begin{equation}
 \dfrac{\partial\vec{X}}{\partial \xi_1} \times
\dfrac{\partial\vec{X}}{\partial \xi_2} \neq \vec{0} \quad \forall \ (\xi_1,\xi_2)\in D\subset \mathbb{R}^2, \ \forall \ 
t\in [0,T],
\end{equation}
and the gap is assumed to be small with regard to the dimension of the bound
surfaces. We take into account that the fluid film between the
surfaces is thin by introducing a small non-dimensional parameter 
$\varepsilon$, and setting that 
\begin{equation}
h^\varepsilon(\xi_1,\xi_2,t) = \varepsilon h(\xi_1,\xi_2,t)
\end{equation}
where 
\begin{equation}
h(\xi_1,\xi_2,t) \ge h_0 > 0, \quad \forall \ (\xi_1,\xi_2)\in D\subset \mathbb{R}^2, \ \forall \ 
t\in [0,T].
\end{equation}

\subsection{Construction of the reference domain}
Let us consider 
\begin{equation}
\Omega=D \times [0,1] \label{eq-Omega}
\end{equation}
a domain independent of
$\varepsilon$ and $t$, which is related to $\Omega^{\varepsilon}_t$ by the following change of variable: 

\begin{eqnarray}
    t^\varepsilon&=&t \label{eq-1-cv} \\
    x_i^\varepsilon &=& x_i(\xi_1,\xi_2,t)+\varepsilon \xi_3 h(\xi_1,\xi_2,t)N_i(\xi_1,\xi_2,t) \label{eq-2-cv}
\end{eqnarray}
where $(\xi_1,\xi_2)\in D$ and $\xi_3 \in[0,1]$.

Let us define the basis $\left\{
\vec{a}_1,\vec{a}_2,\vec{a}_3\right\}$ 
\begin{eqnarray}
 \vec{a}_1(\xi_1,\xi_2,t)&=&\dfrac{\partial \vec{X}(\xi_1,\xi_2,t)}{\partial \xi_1} \label{base_a1} \\
 \vec{a}_2(\xi_1,\xi_2,t)&=&\dfrac{\partial \vec{X}(\xi_1,\xi_2,t)}{\partial \xi_2}\\
 \vec{a}_3(\xi_1,\xi_2,t)&=& \vec{N}(\xi_1,\xi_2,t) \label{base_a3}
\end{eqnarray}

In \ref{ApendiceA} we obtain 
\begin{eqnarray}
\dfrac{ \partial x_i^{\varepsilon}}{\partial \xi_j} &=&
a_{ji}+\varepsilon \xi_3 \dfrac{\partial h}{\partial \xi_j} a_{3i}
+\varepsilon \xi_3 h\dfrac{\partial a_{3i}}{\partial \xi_j}, \quad (i=1,2,3; j=1,2) \label{parcial_xiuv}\\
\dfrac{ \partial x_i^{\varepsilon}}{\partial \xi_3} &=& \varepsilon
h
a_{3i}, \quad (i=1,2,3) \label{parcial_xiw}\\
\dfrac{ \partial x_i^{\varepsilon}}{\partial t} &=& \dfrac{\partial
x_i}{\partial t}+\varepsilon \xi_3 \dfrac{\partial h}{\partial
t}a_{3i} +\varepsilon \xi_3 h\dfrac{\partial a_{3i}}{\partial t},
\quad (i=1,2,3) \label{parcial_xit}  \\
\dfrac{ \partial t^{\varepsilon}}{\partial \xi_1} &=&\dfrac{
\partial t^{\varepsilon}}{\partial \xi_2}=\dfrac{ \partial
t^{\varepsilon}}{\partial \xi_3}=0 , \quad (i=1,2,3)
\label{parcial_tuvw}\\
\dfrac{ \partial t^{\varepsilon}}{\partial t} &=&1,
\label{parcial_tt} 
\end{eqnarray}
where $a_{ij} = \vec{a}_i \cdot \vec{e}_j$, ($i,j=1,2,3$), $\left \{ \vec{e}_1, \vec{e}_2, \vec{e}_3 \right \}$ is the canonical basis of $\mathbb{R}^3$, and 
\begin{eqnarray}
&&\left(\dfrac{ \partial \xi_1}{\partial x_1^{\varepsilon}}, \dfrac{
\partial \xi_1}{\partial x_2^{\varepsilon}}, \dfrac{ \partial \xi_1}{\partial x_3^{\varepsilon}} \right)
= \alpha_1 \vec{a}_1 + \beta_1 \vec{a}_2 +\gamma_1 \vec{a}_3
\label{Dx_xi1_base_a}
\\
&&\left(\dfrac{ \partial \xi_2}{\partial x_1^{\varepsilon}}, \dfrac{
\partial \xi_2}{\partial x_2^{\varepsilon}}, \dfrac{ \partial \xi_2}{\partial x_3^{\varepsilon}} \right)
= \alpha_2 \vec{a}_1 + \beta_2 \vec{a}_2 +\gamma_2 \vec{a}_3
\\
&&\left(\dfrac{ \partial \xi_3}{\partial x_1^{\varepsilon}}, \dfrac{
\partial \xi_3}{\partial x_2^{\varepsilon}}, \dfrac{ \partial \xi_3}{\partial x_3^{\varepsilon}} \right)
= \alpha_3 \vec{a}_1 + \beta_3 \vec{a}_2 +\gamma_3 \vec{a}_3\label{Dx_xi3_base_a}
\end{eqnarray}
\begin{eqnarray}
&&\dfrac{
\partial \xi_1}{\partial t^{\varepsilon}} = -(\alpha_1 \vec{a}_1 + \beta_1 \vec{a}_2)\cdot\left(
\dfrac{\partial \vec{X}}{\partial t} + \varepsilon \xi_3 h
\dfrac{\partial \vec{a}_{3}}{\partial t} \right)
\label{parcial_xi1_teps_alfa_beta}\\
&&\dfrac{
\partial \xi_2}{\partial t^{\varepsilon}} = -( \alpha_2 \vec{a}_1 + \beta_2 \vec{a}_2)\cdot
\left(\dfrac{\partial \vec{X}}{\partial t} + \varepsilon \xi_3 h
\dfrac{\partial \vec{a}_{3}}{\partial t}
\right)\label{parcial_xi2_teps_alfa_beta}\\
&&\dfrac{
\partial \xi_3}{\partial t^{\varepsilon}} = - (\alpha_3 \vec{a}_1 +\beta_3 \vec{a}_2) \cdot\left(\dfrac{\partial \vec{X}}{\partial t}
+\varepsilon \xi_3 h   \dfrac{\partial \vec{a}_{3}}{\partial t}
\right) -\dfrac{1}{\varepsilon h} \vec{a}_3 \cdot\dfrac{\partial
\vec{X}}{\partial t}- \dfrac{\xi_3}{ h} \dfrac{\partial h}{\partial
t} \label{parcial_xi3_teps_alfa_beta}\\
&& \dfrac{\partial t}{\partial x^{\varepsilon}_i}= 0 \quad (i=1,2,3) \label{dtdxi_0} \\
&& \dfrac{\partial t}{\partial t^{\varepsilon}} = 1 \label{dtdteps_1}
\end{eqnarray}
where $\alpha_i$, $\beta_i$, $\gamma_i$ ($i=1,2,3$) are given by
\eqref{alfa1}-\eqref{gamma3} in \ref{ApendiceA}. 

Given any function
$F^\varepsilon(x_1^{\varepsilon},x_2^{\varepsilon},x_3^{\varepsilon},t^{\varepsilon})$
defined on $ \Omega^{\varepsilon}_t$, we can define another function
$F(\varepsilon)(\xi_1,\xi_2,\xi_3,t)$ on $\Omega$ using the change
of variable
\begin{equation}
F(\varepsilon)(\xi_1,\xi_2,\xi_3,t)=F^{\varepsilon}(x_1^{\varepsilon},x_2^{\varepsilon},x_3^{\varepsilon},t^{\varepsilon})
\end{equation}
and the relation between its partial derivatives is trivially:
\begin{eqnarray}
 \dfrac{ \partial F^{\varepsilon}}{\partial x^{\varepsilon}_i} &=&
\displaystyle \frac{
\partial F(\varepsilon)}{\partial \xi_1} \displaystyle \frac{ \partial
\xi_1}{\partial x^{\varepsilon}_i}+ \displaystyle \frac{ \partial
F(\varepsilon)}{\partial \xi_2} \displaystyle \frac{ \partial
\xi_2}{\partial x^{\varepsilon}_i}+ \displaystyle \frac{ \partial
F(\varepsilon)}{\partial \xi_3} \displaystyle \frac{ \partial
\xi_3}{\partial x^{\varepsilon}_i} \label{der_F_xi}\\
   \dfrac{ \partial F^{\varepsilon}}{\partial t^{\varepsilon}} &=&
\displaystyle \frac{ \partial F(\varepsilon)}{\partial t} +
\displaystyle \frac{
\partial F(\varepsilon)}{\partial \xi_1} \displaystyle \frac{ \partial
\xi_1}{\partial t^{\varepsilon}}+ \displaystyle \frac{ \partial
F(\varepsilon)}{\partial \xi_2} \displaystyle \frac{ \partial
\xi_2}{\partial t^{\varepsilon}}+ \displaystyle \frac{ \partial
F(\varepsilon)}{\partial \xi_3} \displaystyle \frac{ \partial
\xi_3}{\partial t^{\varepsilon}}\label{der_F_t}
\end{eqnarray}
where $\dfrac{\partial \xi_j}{\partial x_i^{\varepsilon}}$, $\dfrac{ \partial \xi_j}{\partial t^{\varepsilon}}$ are given by \eqref{Dx_xi1_base_a}-\eqref{parcial_xi3_teps_alfa_beta}.

\subsection{Navier-Stokes Equations}
Let us consider an
incompressible newtonian fluid, so we can assume that the fluid motion
is governed by Navier-Stokes equations ($i=1,2,3$):
\begin{eqnarray}
&& \rho_0\left(\dfrac{\partial  u_i^{\varepsilon}}{\partial
t^{\varepsilon}}+ \dfrac{\partial u_i^{\varepsilon}}{\partial
x_j^{\varepsilon}} u_j^{\varepsilon}\right)=-\dfrac{\partial
p^{\varepsilon}}{\partial x_i^{\varepsilon}}+ \mu \left(
\dfrac{\partial^2 u_i^{\varepsilon}}{\partial (x_1^{\varepsilon})^2}
+
\dfrac{\partial^2  u_i^{\varepsilon}}{\partial (x_2^{\varepsilon})^2} +
 \dfrac{\partial^2  u_i^{\varepsilon}}{\partial (x_3^{\varepsilon})^2}\right)+ \rho_0 f_i^{\varepsilon} \label{ec_ns_ij}\\
&&\dfrac{\partial  u_j^{\varepsilon}}{\partial x_j^{\varepsilon}}=0
\label{div_i}
\end{eqnarray}
where repeated indices indicate summation ($j$ takes values from 1 to 3), $\rho_0$ is the fluid density, assumed to be constant,
$\vec{u}^{\varepsilon}=(u_1^{\varepsilon},u_2^{\varepsilon},u_3^{\varepsilon})$
is the fluid velocity, $p^{\varepsilon}$ is the pressure,
$\mu$ is the dynamic viscosity and  $\vec{f}^\varepsilon$ denotes the external density of 
volume forces.

Let us write $\vec{u}^{\varepsilon}$ and $\vec{f}^\varepsilon$ in the new basis \eqref{base_a1}-\eqref{base_a3} (repeated indices $i$ and $k$ indicate summation from 1 to 3):

\begin{eqnarray}
\vec{u}^{\varepsilon} &=& u_i^{\varepsilon}
\vec{e}_i = u_k(\varepsilon)\vec{a}_{k} \\
\vec{f}^\varepsilon &=& f_i^{\varepsilon}\vec{e}_i = 
f_k(\varepsilon)\vec{a}_{k}
\end{eqnarray}
so we have 
\begin{eqnarray}
u_i^{\varepsilon}&=& \left ( u_k(\varepsilon)\vec{a}_{k} \right ) \cdot \vec{e}_i = u_k(\varepsilon) a_{ki} \label{cambio_base_u}\\
f_i^{\varepsilon}&=& \left ( f_k(\varepsilon)\vec{a}_{k} \right ) \cdot \vec{e}_i = f_k(\varepsilon){a}_{ki}
\label{cambio_base_f}
\end{eqnarray}

Taking into account \eqref{cambio_base_u}-\eqref{cambio_base_f}, equations \eqref{ec_ns_ij}-\eqref{div_i} yield ($i=1,2,3$):
\begin{eqnarray}
&& \rho_0\left(\dfrac{\partial  (u_k(\varepsilon){a}_{ki})}{\partial
t^{\varepsilon}}+ \dfrac{\partial
(u_k(\varepsilon){a}_{ki})}{\partial x_j^{\varepsilon}}
(u_k(\varepsilon){a}_{kj})\right)=-\dfrac{\partial
p(\varepsilon)}{\partial x_i^{\varepsilon}}\nonumber\\
&&{}\hspace*{+0.5cm}+ \mu \left( \dfrac{\partial^2
(u_k(\varepsilon){a}_{ki})}{\partial (x_1^{\varepsilon})^2} +
\dfrac{\partial^2 (u_k(\varepsilon){a}_{ki})}{\partial
(x_2^{\varepsilon})^2} +
 \dfrac{\partial^2  (u_k(\varepsilon){a}_{ki})}{\partial (x_3^{\varepsilon})^2}\right)+
 \rho_0 f_k(\varepsilon){a}_{ki} \label{ec_ns_ij_a}\\
&&\dfrac{\partial  (u_k(\varepsilon){a}_{kj})}{\partial
x_j^{\varepsilon}}=0 \label{div_i_a}
\end{eqnarray}

Equations \eqref{ec_ns_ij_a}-\eqref{div_i_a} can be written in the
reference domain $\Omega$, using
\eqref{der_F_xi}-\eqref{der_F_t} and \eqref{Dx_xi1_base_a}-\eqref{parcial_xi3_teps_alfa_beta}, as follows (repeated indices indicates summation from 1 to 3; $i=1,2,3$):
\begin{eqnarray}
&&\dfrac{ \partial u_k(\varepsilon)}{\partial t}  {a}_{ki} +
u_k(\varepsilon)\dfrac{
\partial {a}_{ki}}{\partial t} \nonumber\\
&&\hspace*{+0.5cm}{}+\left({a}_{ki} \dfrac{
\partial u_k(\varepsilon)}{\partial \xi_l} + u_k(\varepsilon)\dfrac{ \partial {a}_{ki}}{\partial \xi_l} \right)\left[ -(\alpha_l \vec{a}_1 + \beta_l \vec{a}_2)\cdot\left(
\dfrac{\partial \vec{X}}{\partial t} + \varepsilon \xi_3 h
\dfrac{\partial \vec{a}_{3}}{\partial t} \right) \right]
\nonumber\\
&&\hspace*{+0.5cm}{}+\left({a}_{ki}\dfrac{
\partial u_k(\varepsilon)}{\partial \xi_3} + u_k(\varepsilon)\dfrac{
\partial {a}_{ki}}{\partial \xi_3} \right)\left(
-\dfrac{1}{\varepsilon h} \vec{a}_3 \cdot\dfrac{\partial
\vec{X}}{\partial t}- \dfrac{\xi_3}{ h} \dfrac{\partial h}{\partial
t}
\right)\nonumber\\
&&\hspace*{+0.5cm}{}+u_k(\varepsilon){a}_{kj} \left({a}_{ki} \dfrac{
\partial u_k(\varepsilon)}{\partial \xi_l} +u_k(\varepsilon)
\dfrac{
\partial {a}_{ki}}{\partial \xi_l} \right) \left(  \alpha_l
{a}_{1j} + \beta_l {a}_{2j}+ \gamma_l {a}_{3j}\right)\nonumber\\
&&\hspace*{+0.5cm}=-\dfrac{1}{\rho_0}\dfrac{
\partial p(\varepsilon)}{\partial \xi_l}\left(\alpha_l
{a}_{1i} + \beta_l {a}_{2i}+ \gamma_l {a}_{3i}\right)\nonumber\\
&&\hspace*{+0.5cm}{} + \nu \left\{ \left[\dfrac{
\partial^2 (u_k(\varepsilon){a}_{ki})}{\partial \xi_l \partial \xi_m}
\left( \alpha_l {a}_{1j} + \beta_l {a}_{2j}+ \gamma_l {a}_{3j}
\right)\right.\right. \nonumber\\
&&\hspace*{+0.5cm}\left.\left.{}+ \dfrac{
\partial (u_k(\varepsilon){a}_{ki})}{\partial \xi_l}\dfrac{\partial}{\partial \xi_m}\left(
\alpha_l {a}_{1j} + \beta_l {a}_{2j}+ \gamma_l {a}_{3j} \right)
\right] \left( \alpha_m {a}_{1j} + \beta_m {a}_{2j}+ \gamma_m
{a}_{3j} \right) \right\} \nonumber\\
&&\hspace*{+0.5cm}+
f_k(\varepsilon){a}_{ki}\label{ec_ns_ij_alfa_beta}\\
&& \left({a}_{kj} \dfrac{
\partial u_k(\varepsilon)}{\partial \xi_l} + u_k(\varepsilon)\dfrac{
\partial {a}_{kj}}{\partial \xi_l}\right) \left(\alpha_l
{a}_{1j} + \beta_l {a}_{2j}+ \gamma_l {a}_{3j}\right) =0
\label{div_i_dr_alfa_beta}
\end{eqnarray}

\subsection{Asymptotic Analysis}

Let us assume that $u_i(\varepsilon)$, $f_i(\varepsilon)$
($i=1,2,3$) and $p(\varepsilon)$ can be developed in powers of
$\varepsilon$, that is:
\begin{eqnarray}
&& u_i(\varepsilon) = u_i^0 + \varepsilon u_i^1 + \varepsilon^2
u_i^2 + \cdots \quad (i=1,2,3) \label{ansatz_1}\\
&&p(\varepsilon) =\varepsilon^{-2} p^{-2} + \varepsilon^{-1} p^{-1}
+p^0 + \varepsilon p^1 + \varepsilon^2 p^2 + \cdots \label{ansatz_2}\\
&& f_i(\varepsilon) = f_i^0 + \varepsilon f_i^1 + \varepsilon^2
f_i^2 + \cdots \quad (i=1,2,3) \label{ansatz_3}
\end{eqnarray}

In making this choice, we follow \cite{BayadaChambat1986}, \cite{Assemien}, and \cite{Cimatti1987}. 

Before substituting $\alpha_i, \beta_i, \gamma_i\ (i=1,2,3)$ in \eqref{ec_ns_ij_alfa_beta}-\eqref{div_i_dr_alfa_beta}, we must develop \eqref{alfa1}-\eqref{beta312} in powers of $\varepsilon$. It is easy to check that 

\begin{eqnarray}
\alpha_i
&=&
\alpha_i^0+\varepsilon \xi_3 h \alpha_i^1 +\varepsilon^2 \xi_3^2h^2
\alpha_i^2+\cdots,
\quad (i=1,2) \label{alfaides}\\
\alpha_3 &=&\dfrac{ \xi_3 }{ h}(\alpha_3^0+\varepsilon \xi_3 h \alpha_3^1
+\varepsilon^2 \xi_3^2h^2 \alpha_3^2+\cdots), \label{alfa3des}\\
 \beta_i &=&\beta_i^0+\varepsilon \xi_3 h \beta_i^1 +\varepsilon^2 \xi_3^2h^2
\beta_i^2+\cdots, \quad (i=1,2) \label{betaides}\\
\beta_3&=&\dfrac{ \xi_3 }{ h}(\beta_3^0+\varepsilon \xi_3 h \beta_3^1
+\varepsilon^2 \xi_3^2h^2 \beta_3^2+\cdots), \label{beta3des} \\
\gamma_3&=&\dfrac{1}{\varepsilon h}, \quad \gamma_1 = \gamma_2 = 0, \label{gammades} 
\end{eqnarray}
where 
\begin{eqnarray}
&&\alpha_1^0=\dfrac{\|\vec{a}_2\|^2}{A^0}=\dfrac{G}{EG-F^2}\label{alfa10}\\
&&\alpha_1^1= \dfrac{\vec{a}_2\cdot \dfrac{\partial
\vec{a}_{3}}{\partial \xi_2}- \alpha_1^0 A^1}{ A^0}=- \dfrac{g +
\alpha_1^0 A^1}{ A^0}\label{alfa11}
\\
&&\alpha_1^n= -\dfrac{\alpha_1^{n-2} A^2+\alpha_1^{n-1} A^1}{  A^0},
\quad n\geq 2
\label{alfa1n}\\
&&\alpha_2^0=\beta_1^0=-\dfrac{ \vec{a}_2\cdot \vec{a}_{1}}{ A^0} =-\dfrac{ F}{ A^0}  \label{alfa20_beta10}\\
&&\alpha_2^1=\beta_1^1 =- \dfrac{\vec{a}_2\cdot \dfrac{\partial
\vec{a}_{3}}{\partial \xi_1} +\alpha_2^0 A^1}{ A^0}=
\dfrac{f- \alpha_2^0 A^1}{ A^0} \label{alfa21_beta11}\\
&&\alpha_2^n = \beta_1^n= -\dfrac{ \alpha_2^{n-2} A^2 +
\alpha_2^{n-1} A^1}{ A^0},
\quad n\geq 2  \label{alfa2n_beta1n}\\
&&\alpha_3^0=\dfrac{\dfrac{\partial h}{\partial \xi_2} \vec{a}_{1}
\cdot \vec{a}_2 - \dfrac{\partial h}{\partial \xi_1} \|\vec{a}_2\|^2
}{A^0} =\dfrac{\dfrac{\partial h}{\partial \xi_2} F -
\dfrac{\partial h}{\partial \xi_1}
G }{A^0}     \label{alfa30}\\
&&\alpha_3^1= \dfrac{\vec{a}_2\cdot\left[\dfrac{\partial h}{\partial
\xi_2}
  \dfrac{\partial
\vec{a}_{3}}{\partial \xi_1} - \dfrac{\partial h}{\partial \xi_1}
\dfrac{\partial \vec{a}_{3}}{\partial \xi_2}\right]-\alpha_3 ^0 A^1
} { A^0} = \dfrac{-\dfrac{\partial h}{\partial \xi_2} f +
\dfrac{\partial h}{\partial \xi_1} g-\alpha_3 ^0 A^1 } { A^0}
\label{alfa31}\\
 && \alpha_3^n=-\dfrac{\alpha_3^{n-2} A^2 + \alpha_3^{n-1} A^1} {A^0},  \quad n\geq 2 \label{alfa3n}
\end{eqnarray}
\begin{eqnarray}
&&\beta_2^0=\dfrac{ \|\vec{a}_{1}\|^2}{ A^0} =\dfrac{ E}{ A^0}  \label{beta20}\\
&&\beta_2^1 = \dfrac{\vec{a}_1\cdot \dfrac{\partial
\vec{a}_{3}}{\partial \xi_1} -\beta_2^0 A^1}{ A^0} =-
\dfrac{e +\beta_2^0 A^1}{ A^0}\label{beta21}\\
&&\beta_2^n =-\dfrac{ \beta_2^{n-2} A^2 + \beta_2^{n-1} A^1}{ A^0},  \quad n\geq 2   \label{beta22}\\
&&\beta_3^0=\dfrac{\dfrac{\partial h}{\partial \xi_1} \vec{a}_{1}
\cdot \vec{a}_2 - \dfrac{\partial h}{\partial \xi_2}\|\vec{a}_1\|^2
}{A^0} =\dfrac{\dfrac{\partial h}{\partial \xi_1} F -
\dfrac{\partial h}{\partial \xi_2}E }{A^0} \label{beta30}\\
&&\beta_3^1= \dfrac{\dfrac{\partial h}{\partial \xi_1} \left(
\vec{a}_1 \cdot\dfrac{\partial \vec{a}_{3}}{\partial \xi_2}\right) -
\dfrac{\partial h}{\partial \xi_2} \left(  \vec{a}_1
\cdot\dfrac{\partial \vec{a}_{3}}{\partial \xi_1}\right)-\beta_3 ^0
A^1 } { A^0} = \dfrac{-\dfrac{\partial h}{\partial \xi_1} f +
\dfrac{\partial h}{\partial \xi_2} e-\beta_3 ^0 A^1 } { A^0}
\label{beta31}\\
 && \beta_3^n=-\dfrac{\beta_3^{n-2} A^2 + \beta_3^{n-1} A^1} {A^0},  \quad n\geq 2
 \label{beta3n}
\end{eqnarray}

The substitution of the developments \eqref{ansatz_1}-\eqref{ansatz_3} and \eqref{alfaides}-\eqref{beta3n} in \eqref{ec_ns_ij_alfa_beta}-\eqref{div_i_dr_alfa_beta}, and the identification of the terms multiplied by the same power
of $\varepsilon$, lead to a series of equations that will allow
us to determine $\vec{\bf{u}}^0$, $p^{-2}$, etc. 

In what follows, we will use the standard summation convention that repeated indices indicate summation from 1 to 3, unless we indicate otherwise.

In this way, we first identify the terms multiplied by $\varepsilon^{-3}$:
\begin{equation}
-\dfrac{1}{\rho_0}\dfrac{
\partial p^{-2}}{\partial \xi_3} \dfrac{1}{h}
{a}_{3i}=0 \quad (i=1,2,3) 
\end{equation}
so we have
\begin{eqnarray}
\dfrac{
\partial p^{-2}}{\partial \xi_3} =0\label{dp-2_dxi30}
\end{eqnarray}

As a second step,  we identify the terms multiplied by
$\varepsilon^{-2}$.
 Multiplying  by $\vec{a}_i,
\ (i=1,2,3)$, we obtain:
\begin{eqnarray}
&& \dfrac{\mu}{ h^2}\left( E\dfrac{
\partial^2 u_1^0}{\partial \xi_3^2}+ F\dfrac{
\partial^2 u_2^0}{\partial \xi_3^2}\right) = \dfrac{
\partial  p^{-2} }{\partial \xi_1} \label{u10u20p-2a}\\
&& \dfrac{\mu}{ h^2}\left(F \dfrac{
\partial^2 u_1^0}{\partial \xi_3^2} + G \dfrac{
\partial^2 u_2^0}{\partial \xi_3^2}\right) = \dfrac{
\partial  p^{-2} }{\partial \xi_2}\label{u10u20p-2b}\\
&&\dfrac{
\partial   p^{-1}}{\partial \xi_3} \dfrac{1}{ h}=  \dfrac{\mu}{ h^2}\dfrac{
\partial^2 u_3^0}{\partial \xi_3^2} \label{eps-2_ecns3}
\end{eqnarray}

The terms multiplied by $\varepsilon^{-1}$ in \eqref{div_i_dr_alfa_beta} are:
\begin{equation}
\dfrac{
\partial u_k^0}{\partial \xi_3}  {a}_{kj} \dfrac{1}{h} {a}_{3j}= \dfrac{
\partial u_3^0}{\partial \xi_3}   \dfrac{1}{h} =0 \label{du30_xi3_0}
\end{equation}
and using this equality in \eqref{eps-2_ecns3}, we deduce:
\begin{equation}
\dfrac{
\partial   p^{-1}}{\partial \xi_3} =0\label{dp-1_dxi30}
\end{equation}

From the terms multiplied by
$\varepsilon^{-1}$ in \eqref{ec_ns_ij_alfa_beta} we obtain:
\begin{eqnarray}
&&\dfrac{\rho_0 A_0}{ h} \left( E \dfrac{
\partial u_1^0}{\partial \xi_3}+F \dfrac{
\partial u_2^0}{\partial \xi_3} \right) \left(
u_3^0 - \vec{a}_3 \cdot\dfrac{\partial
\vec{X}}{\partial t} \right)\nonumber \\
&&\hspace*{+0.5cm}{}+\xi_3 h\left[\dfrac{
\partial p^{-2}}{\partial \xi_1}   \left(
G e -f F  \right) +\dfrac{
\partial p^{-2}}{\partial \xi_2} \left(f E -e F
\right) \right] + \dfrac{\partial p^{-1}}{\partial \xi_1} A^0  \nonumber \\
&&\hspace*{+0.5cm}= \mu \dfrac{A^0}{ h^2}\left( \dfrac{
\partial^2 u_1^1 }{\partial \xi_3^2} E+ \dfrac{
\partial^2 u_2^1 }{\partial \xi_3^2} F\right)+ \dfrac{\mu A^1}{ h}\left(\dfrac{
\partial u_1^0 }{\partial \xi_3} E + \dfrac{
\partial u_2^0 }{\partial \xi_3} F\right)\label{ns1_eps-1}\\
&&\dfrac{\rho_0 A_0}{ h} \left( F \dfrac{
\partial u_1^0}{\partial \xi_3} + G \dfrac{
\partial u_2^0}{\partial \xi_3}\right) \left(u_3^0- \vec{a}_3
\cdot\dfrac{\partial
\vec{X}}{\partial t} \right)\nonumber \\
&&\hspace*{+0.5cm}{}+\xi_3 h\left[\dfrac{
\partial p^{-2}}{\partial \xi_1}   \left( f G - g F \right)+\dfrac{
\partial p^{-2}}{\partial \xi_2} \left(
 E  g  -f F\right) \right] +
\dfrac{\partial p^{-1}}{\partial \xi_2}A^0\nonumber \\
 &&\hspace*{+0.5cm}  = \mu \dfrac{A^0}{
h^2}\left( \dfrac{
\partial^2 u_1^1 }{\partial \xi_3^2}  F +  \dfrac{
\partial^2 u_2^1 }{\partial \xi_3^2}  G \right)+ \dfrac{\mu A^1}{ h}\left(\dfrac{
\partial u_1^0 }{\partial \xi_3} F + \dfrac{
\partial u_2^0 }{\partial \xi_3} G \right)\label{ns2_eps-1}\\
&&\dfrac{
\partial p^0}{\partial \xi_3}h   = \mu  \dfrac{
\partial^2 u_3^1 }{\partial \xi_3^2} \label{dp0_dxi3}
\end{eqnarray}

Finally, the term of order 0 in \eqref{div_i_dr_alfa_beta} is:
\begin{eqnarray}
&&\dfrac{1}{ h} \dfrac{
\partial u_3^1}{\partial \xi_3} =-\dfrac{\hat{A}_1^0}{A^0}u_1^0-\dfrac{\hat{A}_2^0}{A^0}u_2^0-\dfrac{A^1}{A^0}u_3^0\nonumber\\
&&{}- \dfrac{
\partial u_1^0}{\partial \xi_1}  - \dfrac{
\partial u_2^0}{\partial \xi_2}  +\dfrac{ \xi_3 }{ h } \left(\dfrac{
\partial u_1^0}{\partial \xi_3} \
\dfrac{\partial h}{\partial \xi_1} + \dfrac{
\partial u_2^0}{\partial \xi_3}\dfrac{\partial h}{\partial \xi_2} \right)
\label{ec_div_eps0}
\end{eqnarray}
where
\begin{eqnarray}\hat{A}_i^0 &=& \|\vec{a}_2 \|^2\left(\vec{a}_1
\cdot \dfrac{
\partial \vec{a}_{i}}{\partial \xi_1} \right)- ( \vec{a}_1 \cdot \vec{a}_{2})\left(\vec{a}_2 \cdot \dfrac{
\partial \vec{a}_{i}}{\partial \xi_1} +\vec{a}_1 \cdot \dfrac{
\partial \vec{a}_{i}}{\partial \xi_2}\right)+\|\vec{a}_1\|^2 \left(\vec{a}_2 \cdot \dfrac{
\partial \vec{a}_{i}}{\partial \xi_2} \right)\nonumber\\
&=& \dfrac{1}{2}G\dfrac{\partial E}{\partial
\xi_i}-\dfrac{1}{2}\dfrac{\partial F^2}{\partial \xi_i}+
\dfrac{1}{2}E\dfrac{\partial G}{\partial \xi_i} =
\dfrac{1}{2}\dfrac{\partial(EG- F^2)}{\partial \xi_i}= \dfrac{1}{2}\dfrac{\partial A^0}{\partial \xi_i}, \quad (i=1,2). \label{Bi}
\end{eqnarray}

\section{A new generalized lubrication model} \label{sec-lubrication}

Reynolds wrote, in 1886, a seminal work on lubrication theory (see \cite{Reynolds}), where he introduced heuristically the Reynolds equation. This two-dimensional equation describing the stationary flow of a thin layer of fluid is considered to be the key element for modelling lubrication phenomena. Since then, we can find numerous works in which more general physical models have been considered.

Most models dedicated to the study of thin film flow, specially in lubrication, are derived from the Stokes equation. These first works were focused on stationary models in which the gap and the boundary conditions were fixed with respect to time (see \cite{BayadaChambat1986}, \cite{Cimatti1983}, \cite{Dridi}). These assumptions were considered no longer valid in some devices, so variation with respect to time of the domain was introduced (see   \cite{BayadaChambat1999}). In the same way, 
in some cases, the inertial effects can not be ignored (see \cite{Chaomleffel}), so the studies using Navier-Stokes equation, as ours, turned out to be relevant (see \cite{Assemien}, for example). It was in 1959, in \cite{Elrod}, when full Navier-Stokes equations were used firstly.  Various boundary conditions for the velocity of the surfaces (see \cite{Fantino}), and other types of generalizations have also been studied (see \cite{Chipot1986}, \cite{Chipot1988}, \cite{Cimatti1986} or \cite{Oden1985}).

In this work, as we have stated in the previous section, we will use Navier-Stokes equations to derive a new generalized lubrication model. We are considering a three-dimensional thin domain, that varies with time, whose mean surface can be chosen without any restriction (in particular, neither the lower boundary surface, nor the upper boundary surface, need to be flat). With respect to boundary conditions, we assume that the fluid slips at the lower surface $(\xi_3=0)$, and at the upper surface $(\xi_3=1)$, but there is continuity in the normal direction, so the tangential velocities at the lower and upper surfaces are known, and the normal velocity of each of them must match the fluid velocity.

\begin{eqnarray}
u_k^{\varepsilon} \vec{e}_k = u_k(\varepsilon)\vec{a}_{k}&=&
V_1(\varepsilon) \vec{a}_1+V_2(\varepsilon) \vec{a}_2+\left( \dfrac{\partial
\vec{X}}{\partial t}
\cdot \vec{a}_3\right)\vec{a}_3 \ \textrm{on }\xi_3=0 \label{cc_xi3_0}\\
u_k^{\varepsilon} \vec{e}_k =u_k(\varepsilon)\vec{a}_{k}&=&
W_1(\varepsilon) \vec{a}_1+ W_2(\varepsilon) \vec{a}_2+\left( \dfrac{\partial (\vec{X}+
\varepsilon h\vec{a}_3)}{\partial t} \cdot \vec{a}_3\right)\vec{a}_3
\ \textrm{on }\xi_3=1 \label{cc_xi3_1}
\end{eqnarray}
where $V_1 \vec{a}_1+V_2 \vec{a}_2$ is the tangential velocity
at the lower surface and $W_1\vec{a}_1+W_2\vec{a}_2$ is the
tangential velocity at the upper surface. So we have,
\begin{eqnarray}
u_k(\varepsilon)&=& V_k(\varepsilon) \quad (k=1,2) \quad  \textrm{on }\xi_3=0 \label{cc_uk_xi3_0}\\
u_3(\varepsilon) &=&  \dfrac{\partial \vec{X}}{\partial t}
\cdot \vec{a}_3\quad \textrm{on }\xi_3=0 \label{cc_u3_xi3_0}\\
u_k(\varepsilon)&=& W_k(\varepsilon) \quad (k=1,2)
\quad \textrm{on }\xi_3=1 \label{cc_uk_xi3_1}\\
u_3(\varepsilon)&=& \dfrac{\partial (\vec{X}+\varepsilon
h\vec{a}_3)}{\partial t} \cdot \vec{a}_3 \quad \textrm{on }\xi_3=1
\label{cc_u3_xi3_1}
\end{eqnarray}

If we assume, in the same way as in \eqref{ansatz_1}-\eqref{ansatz_3}, that
\begin{eqnarray}
&& V_i(\varepsilon) = V_i^0 + \varepsilon V_i^1 + \varepsilon^2
V_i^2 + \cdots \quad (i=1,2) \label{ansatzcc0} \\
&& W_i(\varepsilon) = W_i^0 + \varepsilon W_i^1 + \varepsilon^2
W_i^2 + \cdots \quad (i=1,2) \label{ansatzcc}
\end{eqnarray}
we yield from \eqref{cc_xi3_0}-\eqref{cc_xi3_1}:
\begin{eqnarray}
u_k^l&=& V_k^{l} \quad (k=1,2, \quad l=0,1,2,\dots) \quad \textrm{on }\xi_3=0 \label{cc_ukl_xi3_0}\\
u_k^l&=& W_k^{l} \quad (k=1,2, \quad l=0,1,2,\dots) \quad \textrm{on }\xi_3=1 \label{cc_ukl_xi3_1}\\
u_3^0 &=&  \dfrac{\partial \vec{X}}{\partial t}
\cdot \vec{a}_3\quad \textrm{on }\xi_3=0 \label{cc_u30_xi3_0}\\
u_3^l &=&0 \quad (l=1,2,\dots) \quad \textrm{on }\xi_3=0 \label{cc_u3l_xi3_0}\\
u_3^0 &=&  \dfrac{\partial \vec{X}}{\partial t}\cdot \vec{a}_3\quad \textrm{on }\xi_3=1 \label{cc_u30_xi3_1}\\
u_3^1&=& \dfrac{\partial ( h\vec{a}_3)}{\partial t} \cdot \vec{a}_3=
\dfrac{\partial h}{\partial t} \quad \textrm{on }\xi_3=1
\label{cc_u31_xi3_1}\\
u_3^l &=&0 \quad (l=2,3, \dots) \quad \textrm{on }\xi_3=1
\label{cc_u3l_xi3_1}
\end{eqnarray}

From \eqref{u10u20p-2a}-\eqref{u10u20p-2b} we can deduce:
\begin{eqnarray}
&&\hspace*{-0.2cm}\dfrac{
\partial^2 u_2^0}{\partial (\xi_3)^2} =\dfrac{h^2}{\mu A^0}\left( E \dfrac{
\partial p^{-2}}{\partial \xi_2}- F \dfrac{
\partial p^{-2}}{\partial \xi_1}\right)\\
&&\hspace*{-0.2cm} \dfrac{
\partial^2 u_1^0}{\partial (\xi_3)^2} =\dfrac{h^2}{\mu A^0}\left(G \dfrac{
\partial p^{-2}}{\partial \xi_1}- F \dfrac{
\partial p^{-2}}{\partial \xi_2}\right)
\end{eqnarray}

As $p^{-2}$ does not depend on $\xi_3$ (see \eqref{dp-2_dxi30}), we can
integrate the previous equations in $\xi_3$ and impose
\eqref{cc_ukl_xi3_0}-\eqref{cc_ukl_xi3_1}
\begin{eqnarray}
&&\hspace*{-0.2cm} u_1^0 =\dfrac{h^2 (\xi_3^2-\xi_3)}{2\mu
A^0}\left( G \dfrac{
\partial p^{-2}}{\partial \xi_1}- F \dfrac{
\partial p^{-2}}{\partial \xi_2}\right)+\xi_3(W_1^{0}-
V_1^{0})+ V_1^{0}\label{u1_0_lub}\\
&&\hspace*{-0.2cm}
 u_2^0=\dfrac{h^2 (\xi_3^2-\xi_3)}{2\mu A^0}\left(E\dfrac{
\partial p^{-2}}{\partial \xi_2}- F \dfrac{
\partial p^{-2}}{\partial \xi_1}\right)+\xi_3 (W_2^{0}-V_2^{0})+
V_2^{0} \label{u2_0_lub}
\end{eqnarray}

From \eqref{du30_xi3_0}, \eqref{cc_u30_xi3_0} and
\eqref{cc_u30_xi3_1} we know:
\begin{equation}
u_3^0 =  \dfrac{\partial \vec{X}}{\partial t} \cdot \vec{a}_3
\label{u30_lub}
\end{equation}

Now, we yield the following equation by substituting $u^0_i$ ($i=1,2,3$) into equation \eqref{ec_div_eps0} by their
expressions \eqref{u1_0_lub}-\eqref{u30_lub},  integrating over $\xi_3$ from
0 to 1, and evaluating by using \eqref{cc_u3l_xi3_0} and
\eqref{cc_u31_xi3_1}:
\begin{eqnarray}
&&\hspace*{-0.5cm}\dfrac{
\partial }{\partial \xi_1} \left[ \dfrac{h^3 }{ A^0}\left(G \dfrac{
\partial p^{-2}}{\partial \xi_1}- F \dfrac{
\partial p^{-2}}{\partial \xi_2}\right)\right]
+ \dfrac{
\partial }{\partial \xi_2} \left[\dfrac{h^3 }{
A^0}  \left( E \dfrac{
\partial p^{-2}}{\partial \xi_2}- F \dfrac{
\partial p^{-2}}{\partial \xi_1}\right)\right]\nonumber\\
 &&{}=12\mu \dfrac{\partial h}{\partial t} + 12\mu\dfrac{h A^1}{A^0} \left(\dfrac{\partial \vec{X}}{\partial t}
\cdot
\vec{a}_3\right)\nonumber\\
&&{}- \dfrac{h^3 \hat{A}_1^0}{ (A^0)^2} \left( G \dfrac{
\partial p^{-2}}{\partial \xi_1}- F \dfrac{
\partial p^{-2}}{\partial \xi_2}\right)-\dfrac{h^3  \hat{A}_2^0}{ (A^0)^2} \left( E \dfrac{
\partial p^{-2}}{\partial \xi_2}- F \dfrac{
\partial p^{-2}}{\partial \xi_1}\right)  \nonumber\\
&&{}+6\mu\dfrac{h \hat{A}_1^0}{A^0}(W_1^{0}+ V_1^{0}) -
6\mu\dfrac{\partial h}{\partial \xi_1}  (W_1^{0}- V_1^{0}) + 6\mu
h\dfrac{
\partial }{\partial \xi_1}(W_1^{0} +
V_1^{0})\nonumber\\
&&{}+ 6\mu\dfrac{h \hat{A}_2^0}{A^0}  (W_2^{0} + V_2^{0})  - 6\mu
\dfrac{\partial h}{\partial \xi_2} (W_2^{0}-V_2^{0}) + 6\mu h
\dfrac{
\partial }{\partial \xi_2} (W_2^{0} + V_2^{0})
\end{eqnarray}

If we denote by 
\begin{eqnarray}
&&\textrm{div}(f_1,f_2)= \dfrac{
\partial f_1}{\partial \xi_1} + \dfrac{
\partial f_2}{\partial \xi_2}\label{div_def} \\
&&\nabla f= \left(  \dfrac{
\partial f}{\partial \xi_1},  \dfrac{
\partial f}{\partial \xi_2}\right) \label{grad_def}\\
&&\vec{V}^{0}=(V_1^{0},V_2^{0}), \quad
\vec{W}^{0}=(W_1^{0},W_2^{0}), \quad
(\hat{A}_1^0,\hat{A}_2^0)=\dfrac{1}{2}\nabla A^0 \label{vectores}\\
&&M=\begin{pmatrix}G &
-F\\
-F &
E\end{pmatrix}\end{eqnarray} and we take into account that
\begin{equation}
\textrm{div}(\vec{\omega})+\dfrac{1}{2A^0}\nabla A^0 \cdot \vec{\omega} = \dfrac{1}{\sqrt{A^0}}\textrm{div}(\sqrt{A^0}\vec{\omega}) \label{div_A0_w}
\end{equation} we arrive at the equation:
\begin{eqnarray}
&&\hspace*{-0.9cm} \dfrac{1}{\sqrt{A^0}} \textrm{div}\left(\dfrac{h^3 }{ \sqrt{A^0}} M \nabla
p^{-2} \right) =12\mu \dfrac{\partial h}{\partial t} +
12\mu\dfrac{h A^1}{A^0} \left(\dfrac{\partial \vec{X}}{\partial t}
\cdot
\vec{a}_3\right)\nonumber\\
&&\hspace*{-0.5cm}{} - 6\mu \nabla h\cdot (\vec{W}^{0}-\vec{V}^{0}) +
\dfrac{6\mu h}{\sqrt{A^0}} \textrm{div} (\sqrt{A^0}(\vec{W}^{0} + \vec{V}^{0})) \label{Reynolds_gen}
\end{eqnarray}that can be considered a generalization of Reynolds equation.

\begin{rmk}
We claim that \eqref{Reynolds_gen} is a new generalized Reynolds equation because, if we consider the classic assumptions to derive Reynolds equations, we re-obtain the classic Reynolds equation from \eqref{Reynolds_gen}. For example, in \cite{Cimatti1983}, \cite{BayadaChambat1986}, \cite{Cimatti1987} and \cite{Assemien} the domain considered is independent of time, $x_3=0$ in \eqref{eq-o-domain}, the  upper surface is fixed ($\vec{W}=\vec{0}$) and the lower surface is moving in the $x_1$-direction with constant velocity ($\vec{V}=(s,0)$). Under these assumptions, we can choose a surface parametrization, $\vec{X}$, such that $E=G=1$ and $F=0$, and then equation \eqref{Reynolds_gen} writes as the classical Reynolds equation:
\begin{equation}
\textrm{div}\left(h^3 \nabla p^{-2} \right) = 6\mu s \frac{\partial h}{\partial \xi_1} \label{eq-Reynolds-clasica} 
\end{equation}
In \cite{BayadaChambat1999} time is taken into account, allowing $h$ to depend on time, and then the term $\dfrac{\partial h}{\partial t}$ appears: 
\begin{equation}
\textrm{div}\left(h^3 \nabla p^{-2} \right) = 12\mu \frac{\partial h}{\partial t} + 6\mu s \frac{\partial h}{\partial \xi_1} \label{eq-Reynolds-clasica-2} 
\end{equation}
\end{rmk}

\begin{rmk}
The matrix $M$ and the coefficients $A^0$, $A^1$, that appear in \eqref{Reynolds_gen}, depend only on the geometry of the surface parametrized by $\vec{X}$. In fact, the 
matrix $\dfrac{1}{A^0}M$ is the inverse of the matrix of the first fundamental form of $\vec{X}$, and the term $\dfrac{A^1}{A^0}=-2K_{m}$ (see \eqref{Km}).
\end{rmk}

\begin{rmk} \label{rmk-b-c-p}
Equation \eqref{Reynolds_gen} must be completed with boundary conditions at $\partial D$, usually the value of $p^{-2}$ at $\partial D$. 
\end{rmk}

\begin{rmk} \label{rmk-eq-h}
Equation \eqref{Reynolds_gen} can be re-scaled, and then $p^\varepsilon$ is approximated by $p^{-2, \varepsilon} = \varepsilon^{-2} p^{-2}$, solution of 
\begin{eqnarray}
&&\hspace*{-0.9cm} \dfrac{1}{\sqrt{A^0}} \textrm{div}\left(\dfrac{(h^{\varepsilon})^3 }{ \sqrt{A^0}} M \nabla
p^{-2, \varepsilon} \right) =12\mu \dfrac{\partial h^{\varepsilon }}{\partial t^{\varepsilon }} +
12\mu\dfrac{h^{\varepsilon } A^1}{A^0} \left(\dfrac{\partial \vec{X}}{\partial t}
\cdot
\vec{a}_3\right)\nonumber\\
&&\hspace*{-0.5cm}{} - 6\mu \nabla h^{\varepsilon } \cdot (\vec{W}^{0}-\vec{V}^{0}) +
\dfrac{6\mu h^{\varepsilon }}{\sqrt{A^0}} \textrm{div} (\sqrt{A^0}(\vec{W}^{0} + \vec{V}^{0})) \label{Reynolds_gen_resc}
\end{eqnarray}
\end{rmk}

\begin{rmk}
We must point out that the expression 
\begin{equation}
    \dfrac{1}{\sqrt{A^0}}\textrm{div}(\sqrt{A^0}\vec{\omega})= \omega^1_{,1}+\omega^2_{,2}
\end{equation}
is exactly the covariant divergence of $\vec{\omega}$, where $\omega^\alpha_{, \beta}$ stands for the covariant devirative of $\omega^\alpha$ with respect to $\xi_\beta$.
\end{rmk}

\section{A new thin fluid layer model} \label{sec-thin-layer}

Thin fluid layer models are widely used for the analysis and numerical simulation of a large number of geophysical phenomena, such as rivers or coastal flows and other hydraulic applications. Saint-Venant firstly derived in his paper \cite{SaintVenant} a shallow water model, since then numerous authors have studied this type of models (see, for example, \cite{Orenga}, \cite{Sundbye}, \cite{Bresch1}-\cite{BreschNoble}, \cite{FerrariSaleri}, \cite{GerbeauPerthame}, \cite{Marche}), on many occasions using  asymptotic analysis techniques to justify them (see \cite{AzeradGuillen},\cite{RodTab1}-\cite{RodTab6}).

With this aim, in this section we will study what happens when, instead of considering that the tangential and normal velocities are known on the upper and lower surfaces, as we have done in \eqref{cc_xi3_0}-\eqref{cc_xi3_1}, we assume that the normal component of the traction on $\xi_3=0$ and on $\xi_3=1$ are known pressures, and that the tangential component of the traction on these surfaces are friction forces depending on the value of the velocities on $\partial D$. Therefore, we assume that

\begin{eqnarray}
&&\vec{T}^{\varepsilon}\cdot
\vec{n}^{\varepsilon}_0 = (\sigma^{\varepsilon} \vec{n}^{\varepsilon}_0)\cdot
\vec{n}^{\varepsilon}_0=-\pi^{\varepsilon}_0 \textrm{ on } \xi_3=0, \label{Tn0_xi3_0}\\
&&\vec{T}^{\varepsilon}\cdot
\vec{n}^{\varepsilon}_1 = (\sigma^{\varepsilon} \vec{n}^{\varepsilon}_1)\cdot
\vec{n}^{\varepsilon}_1=-\pi^{\varepsilon}_1 \textrm{ on } \xi_3=1,
\label{Tn1_xi3_1}
\\
&&\vec{T}^{\varepsilon}\cdot
\vec{a}_i= (\sigma^{\varepsilon}  \vec{n}^{\varepsilon}_0)\cdot
\vec{a}_i=-\vec{f}^{\varepsilon}_{R0}\cdot
\vec{a}_i \textrm{ on } \xi_3=0,\quad (i=1,2) \label{Tai_xi3_0}\\
&&\vec{T}^{\varepsilon}\cdot
\vec{v}_i^{\varepsilon}= (\sigma^{\varepsilon}   \vec{n}^{\varepsilon}_1)\cdot
\vec{v}_i^{\varepsilon}=-\vec{f}^{\varepsilon}_{R1}\cdot
\vec{v}_i^{\varepsilon} \textrm{ on } \xi_3=1,\quad (i=1,2)
\label{Tvi_xi3_1}
\end{eqnarray}
where $\vec{T}^{\varepsilon}$ is the traction vector and $\sigma^{\varepsilon}$ is the stress tensor given by
\begin{equation}
    \sigma^{\varepsilon}_{ij}= -p^{\varepsilon}\delta_{ij}+\mu \left(
\dfrac{\partial u_i^{\varepsilon}}{\partial x^{\varepsilon}_j} +
\dfrac{\partial u_j^{\varepsilon}}{\partial
x^{\varepsilon}_i}\right),  \quad (i,j=1,2,3)\label{sigmaij}
\end{equation}
vectors $\vec{n}^{\varepsilon}_0$, $\vec{n}^{\varepsilon}_1$ are, respectively, the outward unit normal vectors to the lower and the upper surfaces, that is 
\begin{eqnarray}
&&\vec{n}^{\varepsilon}_0=s_0 \vec{a}_3 \label{n0}\\
&&\vec{n}^{\varepsilon}_1=-s_0\dfrac{\vec{v}^{\varepsilon}_3}{\|\vec{v}^{\varepsilon}_3\|}
\label{n1} 
\end{eqnarray}
where
\begin{equation}
    s_0=-1 \textrm{ or } s_0=1 
\end{equation}
is fixed ($\vec{n}^{\varepsilon}_0 = \vec{a}_3$ or $\vec{n}^{\varepsilon}_0 = - \vec{a}_3$, depending on the orientation of the parametrization $\vec{X}$), and 
\begin{eqnarray}
&&\vec{v}^{\varepsilon}_1=\vec{a}_1+\varepsilon
\left(\dfrac{\partial h}{\partial \xi_1} \vec{a}_3 + h
\dfrac{\partial \vec{a}_3}{\partial \xi_1}\right) \label{v1}\\
&&\vec{v}^{\varepsilon}_2=\vec{a}_2+\varepsilon
\left(\dfrac{\partial h}{\partial \xi_2} \vec{a}_3 + h
\dfrac{\partial \vec{a}_3}{\partial \xi_2}\right) \label{v2}\\
&&\vec{v}^{\varepsilon}_3=\vec{v}^{\varepsilon}_1 \times \vec{v}^{\varepsilon}_2 \label{v3_p}
\end{eqnarray}

From the identities \eqref{v1}-\eqref{v3_p}, we also have the equalities:
\begin{eqnarray}
 &&\vec{v}^{\varepsilon}_3=\vec{a}_1 \times \vec{a}_2 +\varepsilon
\left[\dfrac{\partial h}{\partial \xi_2} (\vec{a}_1 \times
\vec{a}_3) + h \left(\vec{a}_1 \times \dfrac{\partial
\vec{a}_3}{\partial \xi_2}\right) + \dfrac{\partial h}{\partial
\xi_1} (\vec{a}_3\times \vec{a}_2) + h \left(\dfrac{\partial
\vec{a}_3}{\partial \xi_1} \times
\vec{a}_2\right)\right]\nonumber\\
&&\hspace*{0.5cm}+\varepsilon^2 \left[ \left(\dfrac{\partial h}{\partial \xi_1}
\vec{a}_3 + h \dfrac{\partial \vec{a}_3}{\partial
\xi_1}\right)\times\left(\dfrac{\partial h}{\partial \xi_2}
\vec{a}_3 + h \dfrac{\partial \vec{a}_3}{\partial \xi_2}\right)
\right]\label{v3}\\
 &&\|\vec{v}^{\varepsilon}_3\|=\|\vec{a}_1 \times
\vec{a}_2\| + \varepsilon h  \left[ \vec{a}_3 \cdot \left(\vec{a}_1
\times \dfrac{\partial \vec{a}_3}{\partial \xi_2}\right) +\vec{a}_3
\cdot \left(\dfrac{\partial \vec{a}_3}{\partial \xi_1} \times
\vec{a}_2\right)\right] + O(\varepsilon^2)\label{mod_v3}
\end{eqnarray}

Typically, the friction force is of the form
\begin{equation}
    \vec{f}^{\varepsilon}_{R\alpha} = \rho_0 C_R^\varepsilon \|
\vec{{u}}^\varepsilon \| \vec{{u}}^\varepsilon \textrm{ on } \xi_3=\alpha, \quad (\alpha=0,1)
\end{equation}
where $C_R^\varepsilon$ is a small constant. Let us assume that it is of order $\varepsilon$, that is, 
\begin{equation}
    C_R^{\varepsilon} =\varepsilon C^1_R
\end{equation}

Now, taking into account \eqref{sigmaij},
\eqref{Dx_xi1_base_a}-\eqref{Dx_xi3_base_a}, we have the following development in powers of $\varepsilon$:
\begin{eqnarray}
\sigma_{ij}(\varepsilon)&=& -\sum_{r=-2}^{\infty} \varepsilon^r p^r
\delta_{ij} \nonumber\\
&+& \mu \sum_{r=0}^{\infty} \varepsilon^r \left[  \left(\dfrac{
\partial u_k^r}{\partial \xi_l}{a}_{ki} +u_k^r \dfrac{
\partial {a}_{ki}}{\partial \xi_l} \right)\dfrac{ \partial
\xi_l}{\partial x^{\varepsilon}_j}  +  \left(\dfrac{
\partial u_k^r}{\partial \xi_l}{a}_{kj} +u_k^r \dfrac{
\partial {a}_{kj}}{\partial \xi_l}\right)\dfrac{ \partial
\xi_l}{\partial x^{\varepsilon}_i} \right] \nonumber\\
&=& -\varepsilon^{-2} p^{-2} \delta_{ij}+ \varepsilon^{-1} \left\{
-p^{-1} \delta_{ij} +\mu \left[ \dfrac{{a}_{3j}}{ h} \dfrac{
\partial u_k^0}{\partial \xi_3}{a}_{ki} + \dfrac{{a}_{3i}}{ h}  \dfrac{
\partial u_k^0}{\partial \xi_3}{a}_{kj} \right]
\right\} \nonumber\\
&-&p^0 \delta_{ij} + \mu \left[  \dfrac{{a}_{3j}}{ h} \dfrac{
\partial u_k^1}{\partial \xi_3}{a}_{ki} + \dfrac{{a}_{3i}}{ h}  \dfrac{
\partial u_k^1}{\partial \xi_3}{a}_{kj} \right.\nonumber\\
&+&\left. \sum_{l=1}^2 \left(\dfrac{
\partial u_k^0}{\partial \xi_l}{a}_{ki} +u_k^0 \dfrac{
\partial {a}_{ki}}{\partial \xi_l} \right) (\alpha_l^0 {a}_{1j} + \beta_l^0 {a}_{2j})
+ \dfrac{ \xi_3 }{ h}\dfrac{
\partial u_k^0}{\partial \xi_3}{a}_{ki} (\alpha_3^0 {a}_{1j} +
\beta_3^0 {a}_{2j})\right.\nonumber\\
&+&\left. \sum_{l=1}^2 \left(\dfrac{
\partial u_k^0}{\partial \xi_l}{a}_{kj} +u_k^0 \dfrac{
\partial {a}_{kj}}{\partial \xi_l} \right) (\alpha_l^0 {a}_{1i} + \beta_l^0 {a}_{2i})
+ \dfrac{ \xi_3 }{ h}\dfrac{
\partial u_k^0}{\partial \xi_3}{a}_{kj} (\alpha_3^0 {a}_{1i} +
\beta_3^0 {a}_{2i})\right]\nonumber 
\end{eqnarray}
\begin{eqnarray}
&+&\varepsilon \left\{- p^1 \delta_{ij} + \mu \left[
\dfrac{{a}_{3j}}{ h} \dfrac{
\partial u_k^2}{\partial \xi_3}{a}_{ki} + \dfrac{{a}_{3i}}{ h}  \dfrac{
\partial u_k^2}{\partial \xi_3}{a}_{kj}\right. \right.\nonumber\\
&+&\left. \left. \sum_{l=1}^2 \left(\dfrac{
\partial u_k^1}{\partial \xi_l}{a}_{ki} +u_k^1 \dfrac{
\partial {a}_{ki}}{\partial \xi_l} \right) (\alpha_l^0 {a}_{1j} + \beta_l^0 {a}_{2j})
+ \dfrac{ \xi_3 }{ h}\dfrac{
\partial u_k^1}{\partial \xi_3}{a}_{ki} (\alpha_3^0 {a}_{1j} +
\beta_3^0 {a}_{2j})\right. \right.\nonumber\\
&+&\left. \left. \sum_{l=1}^2 \left(\dfrac{
\partial u_k^1}{\partial \xi_l}{a}_{kj} +u_k^1 \dfrac{
\partial {a}_{kj}}{\partial \xi_l} \right) (\alpha_l^0 {a}_{1i} + \beta_l^0 {a}_{2i})
+ \dfrac{ \xi_3 }{ h}\dfrac{
\partial u_k^1}{\partial \xi_3}{a}_{kj} (\alpha_3^0 {a}_{1i} +
\beta_3^0 {a}_{2i})\right. \right.\nonumber\\
&+&\left. \left. \xi_3 h\sum_{l=1}^2 \left(\dfrac{
\partial u_k^0}{\partial \xi_l}{a}_{ki} +u_k^0 \dfrac{
\partial {a}_{ki}}{\partial \xi_l} \right) (\alpha_l^1 {a}_{1j} + \beta_l^1 {a}_{2j})
+ \xi_3^2\dfrac{
\partial u_k^0}{\partial \xi_3}{a}_{ki} (\alpha_3^1 {a}_{1j} +
\beta_3^1 {a}_{2j})\right. \right.\nonumber\\
&+&\left. \left. \xi_3 h\sum_{l=1}^2 \left(\dfrac{
\partial u_k^0}{\partial \xi_l}{a}_{kj} +u_k^0 \dfrac{
\partial {a}_{kj}}{\partial \xi_l} \right) (\alpha_l^1 {a}_{1i} + \beta_l^1 {a}_{2i})
+  \xi_3^2\dfrac{
\partial u_k^0}{\partial \xi_3}{a}_{kj} (\alpha_3^1 {a}_{1i} +
\beta_3^1 {a}_{2i}) \right]\right\} \nonumber\\
&+&\cdots\label{sigma_ij_eps}
\end{eqnarray}

If we assume, now, that
\begin{eqnarray}
\pi_0(\varepsilon)&=&\sum_{r=0}^{\infty} \varepsilon^r \pi_0^r
\label{pi0_serie} \\
\pi_1(\varepsilon)&=&\sum_{r=0}^{\infty} \varepsilon^r \pi_1^r
\label{pi1_serie}\\
\vec{f}_{R0}(\varepsilon)&=&\sum_{r=1}^{\infty} \varepsilon^r
\vec{f}_{R0}^r \label{fR0_serie}\\
\vec{f}_{R1}(\varepsilon)&=&\sum_{r=1}^{\infty} \varepsilon^r
\vec{f}_{R1}^r \label{fR1_serie}
\end{eqnarray}
condition \eqref{Tn0_xi3_0} can be written (using
\eqref{pi0_serie}, \eqref{sigma_ij_eps}, \eqref{n0}) as:
\begin{eqnarray}
&&(\sigma_{ij}(\varepsilon) {a}_{3j}) {a}_{3i}=-\varepsilon^{-2} p^{-2}
-\varepsilon^{-1} p^{-1} - p^0 + \mu \dfrac{2}{ h} \dfrac{
\partial u_3^1}{\partial \xi_3}+\varepsilon \left( -p^1  + \mu  \dfrac{2}{ h} \dfrac{
\partial u_3^2}{\partial \xi_3} \right)+ \cdots\nonumber\\
&&\hspace*{+0.5cm}{} =-(\pi_0^0 + \varepsilon \pi_0^1 + \cdots)\textrm{ on } \xi_3=0
\end{eqnarray}
and we can deduce:
\begin{eqnarray}
&& p^{-2}=0 \textrm{ on } \xi_3=0 \label{p^-2_xi3_0}\\
&& p^{-1}=0 \textrm{ on } \xi_3=0 \label{p^-1_xi3_0}\\
&& - p^0 + \mu \dfrac{2}{ h} \dfrac{
\partial u_3^1}{\partial \xi_3} =-\pi_0^0 \textrm{ on } \xi_3=0 \label{p^0_xi3_0}\\
&& -p^1  + \mu  \dfrac{2}{ h} \dfrac{
\partial u_3^2}{\partial \xi_3}  =- \pi_0^1 \textrm{ on } \xi_3=0
\end{eqnarray}

From \eqref{dp-2_dxi30} and \eqref{p^-2_xi3_0} we obtain 
\begin{eqnarray}
&& p^{-2}=0 \label{p^-2_0}
\end{eqnarray}
and, analogously, from \eqref{dp-1_dxi30} and \eqref{p^-1_xi3_0}, we have 
\begin{eqnarray}
&& p^{-1}=0 \label{p^-1_0}
\end{eqnarray}

Substituting $p^{-2}$ into equations \eqref{u10u20p-2a}-\eqref{u10u20p-2b} we yield 
\begin{equation}
\dfrac{
\partial^2 u_1^0}{\partial \xi_3^2} = \dfrac{
\partial^2 u_2^0}{\partial \xi_3^2} =0 \label{dui0_dxi3_2}
\end{equation}

Let us denote, as in section \ref{sec-lubrication}, by $V_1 \vec{a}_1+V_2 \vec{a}_2$ the tangential velocity
to the lower surface, and by $W_1\vec{a}_1+W_2\vec{a}_2$ the
tangential velocity to the upper surface. Thus we have again \eqref{cc_xi3_0}-\eqref{cc_u3_xi3_1}, but now $V_1(\varepsilon)$, $V_2(\varepsilon)$, $W_1(\varepsilon)$, $W_2(\varepsilon)$, $\dfrac{\partial h}{\partial t}$ are unknown, they are not data as in section \ref{sec-lubrication}.

Let us assume the equalities \eqref{ansatzcc0}-\eqref{ansatzcc} once more. Then we re-obtain \eqref{cc_ukl_xi3_0}-\eqref{cc_u3l_xi3_1}. Now, from \eqref{dui0_dxi3_2} and \eqref{du30_xi3_0}, we deduce 
\begin{eqnarray}
 u_i^0 &=& (W_i^{0}-V_i^{0})\xi_3+V_i^{0} \quad (i=1,2)
\label{u_i^0} \\
u_3^0&=& \dfrac{\partial \vec{X}}{\partial t}
\cdot \vec{a}_3 \label{u_3^0}
\end{eqnarray}

Now, if we substitute $u^0_i$ by their
expressions \eqref{u_i^0}-\eqref{u_3^0} into \eqref{ec_div_eps0}, we integrate over $\xi_3$ from
0 to 1 and we evaluate using \eqref{cc_u3l_xi3_0} and
\eqref{cc_u31_xi3_1}, we obtain
\begin{eqnarray}
2 \dfrac{\partial h}{\partial t} -(\vec{W}^{0}-\vec{V}^{0})\cdot \nabla h +\dfrac{h}{
\sqrt{A^0}}\textrm{div}\left(
\sqrt{A^0}(\vec{W}^{0}+\vec{V}^{0})\right) +\dfrac{2hA^1}{A^0}\left(
\dfrac{\partial \vec{X}}{\partial t}
\cdot \vec{a}_3\right)=0
\label{ec_div_eps0d}
\end{eqnarray}

From \eqref{ns1_eps-1}-\eqref{ns2_eps-1}, \eqref{p^-2_0}-\eqref{p^-1_0} and \eqref{u_i^0}-\eqref{u_3^0}, we have 
\begin{eqnarray}
&&   \dfrac{
\partial^2 u_i^1 }{\partial \xi_3^2} =- \dfrac{h A^1}{A^0} (W_i^0-V_i^0) \qquad (i=1,2) \label{d2ui^1_dxi_3^2} \end{eqnarray}
and integrating twice we yield 
\begin{eqnarray}
&&  u_i^1   =-  \dfrac{h A^1}{2 A^0}(W_i^{0}-V_i^{0})  
(\xi_3^2-\xi_3) + (W_i^{1}-V_i^{1}) \xi_3 + V_i^{1} \quad (i=1,2)
\label{u_i^1}
\end{eqnarray}

From  \eqref{ec_div_eps0}, using
\eqref{u_i^0}, \eqref{u_3^0} and \eqref{div_A0_w}, we can derive an expression for $u_3^1$:
\begin{eqnarray}
&&\hspace*{-1cm}
 u_3^1 =\dfrac{ \xi_3^2 }{2}\left[
(\vec{W}^{0}-\vec{V}^{0}) \cdot \nabla h-\dfrac{ h}{
\sqrt{A^0}}\textrm{div}\left(
\sqrt{A^0}(\vec{W}^{0}-\vec{V}^{0})\right)  \right] \nonumber\\
&&\hspace*{-0.5cm}{}  - h\xi_3\left[\dfrac{1}{\sqrt{A^0}}
\textrm{div}(\sqrt{A^0} \vec{V}^{0})   +\dfrac{A^1}{A^0}\left(
\dfrac{\partial \vec{X}}{\partial t} \cdot \vec{a}_3\right)\right]
\label{u_3^1}
\end{eqnarray}
and we can also yield the following expression for $p^0$ from \eqref{dp0_dxi3}, \eqref{u_3^1} and \eqref{p^0_xi3_0}
\begin{eqnarray}
  p^0   &=& \dfrac{\mu \xi_3}{h}  \left[ -\dfrac{h}{
\sqrt{A^0}}\textrm{div}\left(
\sqrt{A^0}(\vec{W}^{0} - \vec{V}^{0})\right)  +
(\vec{W}^{0}-\vec{V}^{0}) \cdot \nabla h
\right]\nonumber\\
&-&\dfrac{ 2\mu}{\sqrt{A^0}} \textrm{div}(\sqrt{A^0} \vec{V}^{0})
-\dfrac{ 2\mu A^1}{A^0}\left( \dfrac{\partial \vec{X}}{\partial t}
\cdot \vec{a}_3\right)+ \pi_0^0\label{p0}
\end{eqnarray}

Boundary condition \eqref{Tn1_xi3_1} can be written (using
 \eqref{n1}) as follows:
\begin{eqnarray}
&&\left(\sigma_{ij}^{\varepsilon}
{v}^{\varepsilon}_{3j}\right)\cdot
{v}^{\varepsilon}_{3i}=-\pi^{\varepsilon}_1\|\vec{v}^{\varepsilon}_3\|^2
\textrm{ on } \xi_3=1 \label{Tn1_xi3_1s}
\end{eqnarray}

We use expressions  \eqref{sigma_ij_eps} and \eqref{pi1_serie} to substitute $\sigma_{ij}^{\varepsilon}$ and $\pi^{\varepsilon}_1$ into the above condition  and we take into account 
 \eqref{v3}, \eqref{mod_v3}, 
\eqref{p^-2_0}, \eqref{p^-1_0}, \eqref{u_i^0}, \eqref{u_3^1} and
\eqref{p0} to simplify. Identifying the terms multiplied by $\varepsilon^0$ we obtain:
\begin{eqnarray}
&& \|\vec{a}_1 \times \vec{a}_2\|\left\{-\dfrac{\mu }{h}  \left[
\dfrac{h}{ \sqrt{A^0}}\textrm{div}\left(
\sqrt{A^0}(\vec{W}^{0} - \vec{V}^{0})\right)  -
(\vec{W}^{0}-\vec{V}^{0}) \cdot \nabla h \right]
- \pi_0^0 \right\}\nonumber\\
   &&{}+\dfrac{2\mu}{ h}\left[(W_2^{0}-V_2^{0})\left(\dfrac{\partial h}{\partial \xi_2} (\vec{a}_1 \times
\vec{a}_3)\cdot\vec{a}_{2}  + h \left(\vec{a}_1 \times
\dfrac{\partial \vec{a}_3}{\partial \xi_2}\right)\cdot\vec{a}_{2}
\right)\right. \nonumber\\
&&\left.+(W_1^{0}-V_1^{0})\left(  \dfrac{\partial h}{\partial
\xi_1} (\vec{a}_3\times \vec{a}_2) \cdot \vec{a}_1 + h
\left(\dfrac{\partial \vec{a}_3}{\partial \xi_1} \times
\vec{a}_2\right)\cdot
\vec{a}_1\right) \right]\nonumber\\
   && =-\pi_1^0 \|\vec{a}_1 \times
\vec{a}_2\| \textrm{ on } \xi_3=1
\end{eqnarray}

Noticing that
\begin{eqnarray}
&&(\vec{a}_1 \times \vec{a}_3)\cdot\vec{a}_{2}  = (\vec{a}_3\times
\vec{a}_2) \cdot \vec{a}_1 = \left(\vec{a}_1 \times
\dfrac{\vec{a}_{1}\times \vec{a}_2}{\| \vec{a}_{1}\times
\vec{a}_2\|} \right)\cdot \vec{a}_2\nonumber\\
&&\hspace*{+0.5cm} =\dfrac{1}{\| \vec{a}_{1}\times \vec{a}_2\|} \left(\vec{a}_1
(\vec{a}_{1}\cdot \vec{a}_2)- \vec{a}_2 \|\vec{a}_1\|^2
 \right)\cdot \vec{a}_2 = \dfrac{-A^0}{\| \vec{a}_{1}\times \vec{a}_2\|} =- \| \vec{a}_{1}\times \vec{a}_2\| \label{a1xa3_a2}\\
&& \left(\vec{a}_1 \times \dfrac{\partial \vec{a}_3}{\partial
\xi_2}\right)\cdot\vec{a}_{2} =( \vec{a}_{2}\times \vec{a}_1 ) \cdot
\dfrac{\partial \vec{a}_3}{\partial \xi_2} =-\| \vec{a}_{2}\times
\vec{a}_1\| \vec{a}_3 \cdot
\dfrac{\partial \vec{a}_3}{\partial \xi_2}=0 \label{a1_x_da3dxi2_a2}\\
&&\left(\dfrac{\partial \vec{a}_3}{\partial \xi_1} \times
\vec{a}_2\right)\cdot \vec{a}_1 = \dfrac{\partial
\vec{a}_3}{\partial \xi_1} \cdot ( \vec{a}_{2}\times \vec{a}_1 ) =-
\dfrac{\partial \vec{a}_3}{\partial \xi_1} \cdot \|
\vec{a}_{2}\times \vec{a}_1\| \vec{a}_3=0 \label{a2_x_da3dxi1_a1}
\end{eqnarray}
we finally derive
\begin{eqnarray}
&&   \dfrac{\mu}{ \sqrt{A^0}}\textrm{div}\left(
\sqrt{A^0}(\vec{W}^{0}-\vec{V}^{0})\right)  +\dfrac{\mu }{h}
(\vec{W}^{0}-\vec{V}^{0}) \cdot \nabla h + \pi_0^0
 =\pi_1^0 \label{relacion_pi00_pi10}
\end{eqnarray}

Boundary conditions \eqref{Tai_xi3_0} on $\xi_3=0$ can be written using  \eqref{alfa10}, \eqref{alfa20_beta10},
 \eqref{beta20}, \eqref{alfa30}, \eqref{beta30}, \eqref{sigma_ij_eps} and  \eqref{fR0_serie}  in this way:
  \begin{eqnarray}
&&\hspace*{-0.5cm}\varepsilon^{-1}\dfrac{\mu}{ h}  \left( E \dfrac{
\partial u_1^0}{\partial \xi_3}  + F \dfrac{
\partial u_2^0}{\partial \xi_3}
\right) + \mu \left[  \dfrac{E}{ h} \dfrac{
\partial u_1^1}{\partial \xi_3} + \dfrac{F}{ h} \dfrac{
\partial u_2^1}{\partial \xi_3}  + \dfrac{
\partial u_3^0}{\partial \xi_1} + e u_1^0 + f u_2^0
\right]\nonumber\\
&&{}+\varepsilon \mu \left[ \dfrac{E}{ h} \dfrac{
\partial u_1^2}{\partial \xi_3} +\dfrac{F}{ h} \dfrac{
\partial u_2^2}{\partial \xi_3}+ \dfrac{
\partial u_3^1}{\partial \xi_1}+ e u_1^1 + f u_2^1
- \dfrac{ \xi_3 }{ h}\dfrac{
\partial u_3^1}{\partial \xi_3}\dfrac{\partial h}{\partial \xi_1} \right.\nonumber\\
&&\left.{}+ \xi_3 h\sum_{l=1}^2 \left(\dfrac{
\partial u_3^0}{\partial \xi_l} +u_k^0 \dfrac{
\partial \vec{a}_{k}}{\partial \xi_l} \cdot \vec{a}_3\right) (\alpha_l^1 E + \beta_l^1
F)
 \right] +\cdots\nonumber\\
 && =-s_0\left(\varepsilon\vec{f}^1_{R0} + \varepsilon^2 \vec{f}^2_{R0} + \cdots\right)\cdot \vec{a}_{1} \textrm{ on }
\xi_3=0 \label{desarrollo_cc_fR01}\end{eqnarray}
  \begin{eqnarray}
&&\hspace*{-0.5cm} \varepsilon^{-1}\dfrac{\mu }{ h} \left(F \dfrac{
\partial u_1^0}{\partial \xi_3} + G\dfrac{
\partial u_2^0}{\partial \xi_3}
\right)+ \mu \left[  \dfrac{F}{ h} \dfrac{
\partial u_1^1}{\partial \xi_3}+ \dfrac{G}{ h} \dfrac{
\partial u_2^1}{\partial \xi_3}  + \dfrac{
\partial u_3^0}{\partial \xi_2} + f u_1^0 + g u_2^0 
\right]\nonumber\\
&&{}+\varepsilon \mu \left[ \dfrac{F}{ h} \dfrac{
\partial u_1^2}{\partial \xi_3}+ \dfrac{G}{ h} \dfrac{
\partial u_2^2}{\partial \xi_3}+\dfrac{
\partial u_3^1}{\partial \xi_2} + f u_1^1 + g u_2^1  
 - \dfrac{ \xi_3 }{ h}\dfrac{
\partial u_3^1}{\partial \xi_3} \dfrac{\partial h}{\partial \xi_2} \right.\nonumber\\
&&\left.{}+  \xi_3 h\sum_{l=1}^2 \left(\dfrac{
\partial u_3^0}{\partial \xi_l} +u_k^0 \dfrac{
\partial \vec{a}_{k}}{\partial \xi_l} \cdot \vec{a}_3\right) (\alpha_l^1  F + \beta_l^1 G)
 \right]+\cdots\nonumber\\
&&{}=-s_0\left(\varepsilon\vec{f}^1_{R0} + \varepsilon^2 \vec{f}^2_{R0}
+ \cdots\right)\cdot \vec{a}_{2} \textrm{ on } \xi_3=0 \label{desarrollo_cc_fR02}
\end{eqnarray}

We yield from the terms multiplied by $\varepsilon^{-1}$ in the equations above and the equality \eqref{u_i^0} that
\begin{eqnarray}
&&u_i^0=W_i^{0}=V_i^{0} \quad (i=1,2) \label{u_i^0s}
\end{eqnarray}

Identifying the terms multiplied by $\varepsilon^0$  in \eqref{desarrollo_cc_fR01}-\eqref{desarrollo_cc_fR02}, and taking into account \eqref{u_3^0}, \eqref{u_i^1} and \eqref{u_i^0s}, we have
\begin{eqnarray}
&& \hspace*{-1.2cm}    W_1^{1}-V_1^{1}  =-h \left[\sum_{l=1}^2\alpha_l^0 \dfrac{
\partial}{\partial \xi_l}\left(  \dfrac{\partial \vec{X}}{\partial t}
\cdot \vec{a}_3\right)  +  \dfrac{B_{11}}{A^0} V_1^0 + \dfrac{B_{12}}{A^0} V_2^0 
\right]  \label{V_1^1_V^0_abrev}
\\
&&\hspace*{-1.2cm}   W_2^{1}-V_2^{1}  =-h  \left[  \sum_{l=1}^2\beta_l^0\dfrac{
\partial}{\partial \xi_l}\left(  \dfrac{\partial \vec{X}}{\partial t}
\cdot \vec{a}_3\right) +\dfrac{B_{21}}{A^0} V_1^0  + \dfrac{B_{22}}{A^0} V_2^0 
\right]  \label{V_2^1_V^0_abrev}
\end{eqnarray}
where $B_{ij} \ (i,j=1,2)$, are given by \eqref{B11}-\eqref{B22}.

Equalities \eqref{u_i^0s}-\eqref{V_2^1_V^0_abrev} allow us to simplify \eqref{ec_div_eps0d}, \eqref{u_i^1}-\eqref{p0} and \eqref{relacion_pi00_pi10}:
\begin{eqnarray}
&&\hspace*{-1cm}\dfrac{\partial h}{\partial t} +\dfrac{h}{
\sqrt{A^0}}\textrm{div}\left(
\sqrt{A^0} \vec{V}^{0}\right) +\dfrac{hA^1}{A^0}\left(
\dfrac{\partial \vec{X}}{\partial t}
\cdot \vec{a}_3\right)=0
\label{ec_div_eps0ds}
\\
&& \hspace*{-1cm} u_1^1   = -h\left[\dfrac{B_{11}}{A^0} V_1^0  +\dfrac{B_{12}}{A^0} V_2^0  + \sum_{l=1}^2\alpha_l^0  \dfrac{
\partial}{\partial \xi_l}\left(  \dfrac{\partial \vec{X}}{\partial t}
\cdot \vec{a}_3\right) \right] \xi_3 + V_1^{1}
\label{u_1^1s}\\
&&\hspace*{-1cm}  u_2^1  = -h \left[ \dfrac{B_{21}}{A^0}V_1^0  + \dfrac{B_{22}}{A^0} V_2^0  + \sum_{l=1}^2\beta_l^0 \dfrac{
\partial}{\partial \xi_l}\left(  \dfrac{\partial \vec{X}}{\partial t}
\cdot \vec{a}_3\right) \right] \xi_3 + V_2^{1} 
\label{u_2^1s}\\
&&\hspace*{-1cm}
 u_3^1 =  - h\xi_3\left[\dfrac{1}{\sqrt{A^0}}
\textrm{div}(\sqrt{A^0} \vec{V}^{0})   +\dfrac{A^1}{A^0}\left(
\dfrac{\partial \vec{X}}{\partial t} \cdot \vec{a}_3\right)\right]
\label{u_3^1s}\\
&&=\xi_3 \dfrac{\partial h}{\partial t}\label{u_3^1_dhdt}\\
&&\hspace*{-1cm}  p^0  =-\dfrac{ 2\mu}{\sqrt{A^0}} \textrm{div}(\sqrt{A^0} \vec{V}^0)
-\dfrac{ 2\mu A^1}{A^0}\left( \dfrac{\partial \vec{X}}{\partial t}
\cdot \vec{a}_3\right)+ \pi_0^0 \label{p0s}\\
&&=\frac{2\mu}{h} \dfrac{\partial h}{\partial t} + \pi_0^0 \label{p0s-v2} \\
&&\hspace*{-1cm}\pi_1^0=\pi_0^0 \label{pi_0^0=pi_1^0}
\end{eqnarray}

Now, we identify the terms multiplied by $\varepsilon$ in \eqref{desarrollo_cc_fR01}-\eqref{desarrollo_cc_fR02} and, considering \eqref{u_1^1s}-\eqref{u_3^1s}, we obtain:
\begin{eqnarray}
&&\hspace*{-0.5cm}  \mu \left[ \dfrac{A^0}{ h} \dfrac{
\partial u_1^2}{\partial \xi_3} + B_{11} V_1^{1} 
 + B_{12}  V_2^{1} 
 \right] = -s_0\left[G(\vec{f}^1_{R0} \cdot \vec{a}_{1})- 
F(\vec{f}^1_{R0}\cdot \vec{a}_{2})\right] \textrm{ on } \xi_3=0 \label{du_1^2_dxi3_xi3_0}\\
&&\hspace*{-0.5cm} \mu \left[ \dfrac{A^0}{ h} \dfrac{
\partial u_2^2}{\partial \xi_3}+ B_{21} 
 V_1^{1}
+ B_{22}  V_2^{1}
  \right]= - s_0[E(\vec{f}^1_{R0}\cdot \vec{a}_{2}) 
  - F  (\vec{f}^1_{R0} \cdot \vec{a}_{1})] \textrm{ on } \xi_3=0\label{du_2^2_dxi3_xi3_0}
\end{eqnarray}

Going back to \eqref{Tn1_xi3_1s}, the terms multiplied by $\varepsilon$ yield 
\begin{eqnarray}
&&\pi_1^1 = p^1 -\dfrac{2\mu}{ h} \dfrac{
\partial u_3^2}{\partial \xi_3}+  \dfrac{2\mu }{ h } \left[ 
 \dfrac{\partial h}{\partial \xi_2}  (W_2^{1}-V_2^{1})+ \dfrac{\partial h}{\partial \xi_1}
 (W_1^{1}-V_1^{1}) \right]\nonumber\\
&&\hspace*{+0.5cm}{}+2\mu \left[\sum_{l=1}^2 \dfrac{
\partial }{\partial \xi_l}\left(  \dfrac{\partial \vec{X}}{\partial t}
\cdot \vec{a}_3\right) \left(\alpha_l^0 \dfrac{\partial h}{\partial
\xi_1} + \beta_l^0 
\dfrac{\partial h}{\partial \xi_2} \right) \right. \nonumber \\
&&\hspace*{+0.5cm} \left.{}+ V_1^0 \left(\dfrac{\partial h}{\partial
\xi_1} \dfrac{B_{11}}{A^0} + \dfrac{\partial h}{\partial
\xi_2} \dfrac{B_{21}}{A^0}\right) +  V_2^0 \left( \dfrac{\partial h}{\partial
\xi_1} \dfrac{B_{12}}{A^0} + \dfrac{\partial h}{\partial
\xi_2} \dfrac{ B_{22} } {A^0} \right) \right] 
\textrm{ on } \xi_3=1 \label{eq-pi-1-1}
\end{eqnarray}
and using \eqref{V_1^1_V^0_abrev}-\eqref{V_2^1_V^0_abrev}, we can simplify \eqref{eq-pi-1-1}, and write 
\begin{eqnarray}
&& p^1=\pi_1^1 +\dfrac{2\mu}{ h} \dfrac{
\partial u_3^2}{\partial \xi_3}
\textrm{ on } \xi_3=1 \label{p1_xi3_1}
\end{eqnarray}

 Taking into account \eqref{alfa10}, \eqref{alfa20_beta10}, \eqref{beta20}, \eqref{v1}-\eqref{mod_v3},  \eqref{sigma_ij_eps}, \eqref{fR1_serie}, \eqref{u_3^0},  \eqref{a1xa3_a2}-\eqref{a2_x_da3dxi1_a1}, \eqref{u_i^0s}, \eqref{V_1^1_V^0_abrev}-\eqref{V_2^1_V^0_abrev} and \eqref{u_1^1s}-\eqref{u_3^1_dhdt}, we can rewrite conditions  \eqref{Tvi_xi3_1}, identify the terms of order zero, 
\begin{eqnarray}
&&  \dfrac{E}{ h} \dfrac{
\partial u_1^1}{\partial \xi_3}  +  \dfrac{F}{ h} \dfrac{
\partial u_2^1}{\partial \xi_3} +\dfrac{
\partial u_3^0}{\partial \xi_1} + e V_1^0 
+f V_2^0 
=0 \textrm{ on } \xi_3=1 \label{f_R1_eps0}\\
&&  \dfrac{F}{ h} \dfrac{
\partial u_1^1}{\partial \xi_3}  +  \dfrac{G}{ h} \dfrac{
\partial u_2^1}{\partial \xi_3} + \dfrac{
\partial u_3^0}{\partial \xi_2} + f V_1^0 
 + g V_2^0 =0 \textrm{ on } \xi_3=1  \label{f_R12_eps0}
\end{eqnarray}
and the first order terms (here, repeated index $k$ indicates sum from 1 to 3),
\begin{eqnarray}
&&\hspace*{-0.5cm}\mu \left[
\dfrac{A^0}{ h} \dfrac{
\partial u_1^2}{\partial \xi_3}  +  G\dfrac{
\partial^2 h}{\partial t \partial \xi_1}  - F\dfrac{
\partial^2 h}{\partial t \partial \xi_2}  + B_{11} V_1^{1} + B_{12}  V_2^{1} - \dfrac{1 }{ h}\dfrac{
\partial h}{\partial t} \left(\dfrac{\partial h}{\partial \xi_2} F -
\dfrac{\partial h}{\partial \xi_1}
G \right)  \right.\nonumber\\
&&\left.{}- A^0\sum_{l=1}^2 \left[\dfrac{
\partial V_1^0}{\partial \xi_l}
+ \sum_{m=1}^2 \alpha_m^0  u_k^0 \left(\dfrac{
\partial \vec{a}_{k}}{\partial \xi_l} \cdot \vec{a}_m\right)\right] \left( \alpha_l^0  \dfrac{\partial h}{\partial
\xi_1}  + \beta_l^0 \dfrac{\partial h}{\partial \xi_2}\right) 
\right.\nonumber\\
&&\left.{} -A^0 \sum_{l=1}^2 \left(\dfrac{\partial h}{\partial \xi_2} \dfrac{
\partial V_2^0}{\partial \xi_l} + \dfrac{\partial h}{\partial \xi_1} \dfrac{
\partial V_1^0}{\partial \xi_l} \right)\alpha_l^0\right] 
\nonumber\\
&&{}+ \dfrac{\mu}{\sqrt{A^0}} \left[-h I \left( B_{11} V_1^0 + B_{12} V_2^0 
\right)   +A^0 u_k^0\sum_{l=1}^2 \alpha_l^0 \left(\dfrac{
\partial \vec{a}_{k}}{\partial \xi_l} \cdot
\vec{\eta}(h)\right)\right] \nonumber\\
&&=s_0 \left( G \vec{f}^1_{R1}  \cdot
 \vec{a}_1  -F \vec{f}^1_{R1}  \cdot
 \vec{a}_2 \right)
 \textrm{ on } \xi_3=1 \label{du_1^2_dxi3_xi3_1}
\end{eqnarray}

\begin{eqnarray}
&&\hspace*{-0.5cm}\mu \left[
 \dfrac{A^0}{ h} \dfrac{
\partial u_2^2}{\partial \xi_3}  -F \dfrac{
\partial^2 h}{\partial t \partial \xi_1} + E \dfrac{
\partial^2 h}{\partial t \partial \xi_2}  + B_{21} V_1^{1}  +  B_{22} V_2^{1} 
- \dfrac{1 }{ h}\dfrac{
\partial h}{\partial t} \left(\dfrac{\partial h}{\partial \xi_1} F -
\dfrac{\partial h}{\partial \xi_2}E \right)\right.\nonumber\\
&&\left.{}-A^0 \sum_{l=1}^2 \left[\dfrac{
\partial V_2^0}{\partial \xi_l}
+ \sum_{m=1}^2 \beta_m^0  u_k^0 \left(\dfrac{
\partial \vec{a}_{k}}{\partial \xi_l} \cdot \vec{a}_m\right)\right] \left( \alpha_l^0  \dfrac{\partial h}{\partial
\xi_1}  + \beta_l^0 \dfrac{\partial h}{\partial \xi_2}\right)\right.\nonumber\\
&&\left.{} -A^0 \sum_{l=1}^2 \left(\dfrac{\partial h}{\partial \xi_2} \dfrac{
\partial V_2^0}{\partial \xi_l} + \dfrac{\partial h}{\partial \xi_1} \dfrac{
\partial V_1^0}{\partial \xi_l} \right)\beta_l^0 \right] 
\nonumber\\
&&{}+\dfrac{ \mu}{\sqrt{A^0}} \left[ -h I \left( B_{21} V_1^0  + B_{22} V_2^0 
 \right) 
 +A^0 u_k^0 \sum_{l=1}^2 \beta_l^0 \left(\dfrac{
\partial \vec{a}_{k}}{\partial \xi_l} \cdot
\vec{\eta}(h)\right)\right]  \nonumber\\
&&=s_0 \left( - F\vec{f}^1_{R1}  \cdot
 \vec{a}_1  +E \vec{f}^1_{R1}  \cdot
 \vec{a}_2 \right)
 \textrm{ on } \xi_3=1 \label{du_2^2_dxi3_xi3_1}
\end{eqnarray}
where $I$ and $\vec{\eta}(h)$ are given by \eqref{I} and \eqref{eta}.

From equalities \eqref{du_1^2_dxi3_xi3_0}-\eqref{du_2^2_dxi3_xi3_0} and \eqref{du_1^2_dxi3_xi3_1}-\eqref{du_2^2_dxi3_xi3_1} we have (again repeated index $k$ indicates sum from 1 to 3):
\begin{eqnarray}
&&\mu \left[
\dfrac{A^0}{ h} \left(\left. \dfrac{
\partial u_1^2}{\partial \xi_3}\right|_{\xi_3=1} -\left. \dfrac{
\partial u_1^2}{\partial \xi_3}\right|_{\xi_3=0} \right)  +  G\dfrac{
\partial^2 h}{\partial t \partial \xi_1}  - F\dfrac{
\partial^2 h}{\partial t \partial \xi_l} 
- \dfrac{1 }{ h}\dfrac{
\partial h}{\partial t} \left(\dfrac{\partial h}{\partial \xi_2} F -
\dfrac{\partial h}{\partial \xi_1}
G \right) \right.\nonumber\\
&&\left.{}-A^0 \sum_{l=1}^2 \left[\dfrac{
\partial V_1^0}{\partial \xi_l}
+ \sum_{m=1}^2 \alpha_m^0  u_k^0 \left(\dfrac{
\partial \vec{a}_{k}}{\partial \xi_l} \cdot \vec{a}_m\right)\right] \left( \alpha_l^0  \dfrac{\partial h}{\partial
\xi_1}  + \beta_l^0 \dfrac{\partial h}{\partial \xi_2}\right) 
\right.\nonumber\\
&&\left.{} - A^0\sum_{l=1}^2 \left(\dfrac{\partial h}{\partial \xi_2} \dfrac{
\partial V_2^0}{\partial \xi_l} + \dfrac{\partial h}{\partial \xi_1} \dfrac{
\partial V_1^0}{\partial \xi_l} \right)\alpha_l^0\right] 
\nonumber\\
&&{}+ \dfrac{\mu}{\sqrt{A^0}} \left[-h I \left(B_{11} V_1^0 + B_{12} V_2^0 
\right)  +A^0 u_k^0\sum_{l=1}^2 \alpha_l^0 \left(\dfrac{
\partial \vec{a}_{k}}{\partial \xi_l} \cdot
\vec{\eta}(h)\right)\right]  \nonumber\\
&&=s_0\left[ G\left(\vec{f}^1_{R1} +\vec{f}^1_{R0} \right) \cdot
 \vec{a}_1 -F \left(\vec{f}^1_{R1} +\vec{f}^1_{R0} \right) \cdot
 \vec{a}_2 \right]\label{dif_du_1^2_dxi3_1_0}
\end{eqnarray}

\begin{eqnarray}
&&\mu \left[ \left.
 \dfrac{A^0}{ h}\left( \dfrac{
\partial u_2^2}{\partial \xi_3}\right|_{\xi_3=1} -\left.
 \dfrac{
\partial u_2^2}{\partial \xi_3}\right|_{\xi_3=0}  \right) -F \dfrac{
\partial^2 h}{\partial t \partial \xi_1} + E \dfrac{
\partial^2 h}{\partial t \partial \xi_2}
- \dfrac{1 }{ h}\dfrac{
\partial h}{\partial t} \left(\dfrac{\partial h}{\partial \xi_1} F -
\dfrac{\partial h}{\partial \xi_2}E \right)\right.\nonumber\\
&&\left.{}- A^0\sum_{l=1}^2 \left[\dfrac{
\partial V_2^0}{\partial \xi_l}
+ \sum_{m=1}^2 \beta_m^0  u_k^0 \left(\dfrac{
\partial \vec{a}_{k}}{\partial \xi_l} \cdot \vec{a}_m\right)\right] \left( \alpha_l^0  \dfrac{\partial h}{\partial
\xi_1}  + \beta_l^0 \dfrac{\partial h}{\partial \xi_2}\right)\right.\nonumber\\
&&\left.{} -A^0 \sum_{l=1}^2 \left(\dfrac{\partial h}{\partial \xi_2} \dfrac{
\partial V_2^0}{\partial \xi_l} + \dfrac{\partial h}{\partial \xi_1} \dfrac{
\partial V_1^0}{\partial \xi_l} \right)\beta_l^0 \right] 
\nonumber\\
&&{}+\dfrac{ \mu}{\sqrt{A^0}} \left[   -h I \left( B_{21} V_1^0  + B_{22} V_2^0 
 \right)   
 +A^0 u_k^0 \sum_{l=1}^2 \beta_l^0 \left(\dfrac{
\partial \vec{a}_{k}}{\partial \xi_l} \cdot
\vec{\eta}(h)\right)\right]  \nonumber\\
&&=s_0\left[  -F \left(\vec{f}^1_{R1} +\vec{f}^1_{R0} \right) \cdot
 \vec{a}_1 + E \left(\vec{f}^1_{R1} +\vec{f}^1_{R0} \right) \cdot
 \vec{a}_2 \right]\label{dif_du_2^2_dxi3_1_0}
\end{eqnarray}

Now, from the terms of order $\varepsilon^0$ in the equation \eqref{ec_ns_ij_alfa_beta}, and following the steps outlined in \ref{ApendiceC}, we obtain the equations below,

\begin{eqnarray}
&&\hspace*{-0.5cm}\dfrac{ \partial V_i^0}{\partial t}  + \sum_{l=1}^2 \left( V_l^0-C^0_l \right) \dfrac{
\partial V_i^0}{\partial \xi_l} +
\sum_{k=1}^2 
\left( R^0_{ik}+\sum_{l=1}^2
H^0_{ilk} V_l^0 \right) V_k^0\nonumber\\
&&=-\dfrac{1}{\rho_0}\left( \alpha_i^0 \dfrac{
\partial \pi_0^0}{\partial \xi_1}  + \beta_i^0 \dfrac{
\partial \pi_0^0}{\partial \xi_2} \right) \nonumber \\ 
&&{} +\nu \left\{ \sum_{m=1}^2 \sum_{l=1}^2 
 \dfrac{
\partial^2 V_i^0 }{\partial \xi_m \partial \xi_l} J^0_{lm}
+ \sum_{k=1}^2 \sum_{l=1}^2  \dfrac{
\partial V_k^0 }{\partial \xi_l}( L^0_{kli} +\psi(h)^0_{ikl} )\right.\nonumber\\
&& \left. {}
+ \sum_{k=1}^2 V_k^0 ({S}_{ik}^0+\chi(h)^0_{ik})  + \hat{\kappa}(h)^0_i  \right\}+ {F}^0_i(h)-
 Q^0_{i3} \left ( \frac{\partial \vec{X}}{\partial t}
\cdot \vec{a}_3\right )  \quad (i=1,2)\label{ec_Vi0-v2} \end{eqnarray}
where the different coefficients are defined in \ref{ApendiceB}.

\begin{rmk}
Equations \eqref{ec_Vi0-v2} and \eqref{ec_div_eps0ds} allow us to determine $h$, $V_1^0$ and $V_2^0$, once the initial and boundary conditions have been set. These equations provide a shallow water model (see \cite{Orenga}, \cite{Sundbye}, \cite{BreschNoble}, \cite{Marche}, \cite{RodTab1}-\cite{RodTab6}). Equation \eqref{ec_div_eps0ds} represents the conservation of mass of the fluid.  If $h$ is known, then \eqref{ec_div_eps0ds} means an additional condition on the velocity $\vec{V}^0$ and, in that case, the pressure $\pi_0^0$  must be an unknown in \eqref{ec_Vi0-v2}.
\end{rmk}

\begin{rmk}
As in \eqref{Reynolds_gen_resc} (see remark \ref{rmk-eq-h}) equations \eqref{ec_Vi0-v2} and \eqref{ec_div_eps0ds} can be re-scaled, to work with $h^\varepsilon$ instead of $h$.
\end{rmk}

\section{Conclusions \label{conc}}

In this paper, starting from the same initial problem, an incompressible viscous fluid moving between two surfaces parametrized by $\vec{X}$ and $\vec{X} + h^\varepsilon \vec{N}$  (see section \ref{sec-domain}), we obtain, using the asymptotic expansion technique, two different models. The first one is yielded in section \ref{sec-lubrication}, assuming that  the fluid velocity is known on the surfaces $\vec{X}$ and $\vec{X} + h^\varepsilon \vec{N}$. The second one is derived in section \ref{sec-thin-layer}, assuming that we know the tractions applied on the surfaces $\vec{X}$ and $\vec{X} + h^\varepsilon \vec{N}$, rather than the fluid velocity, as in section \ref{sec-lubrication}. This simple change gives rise to two different models: a lubrication model in section \ref{sec-lubrication} and a shallow water model in section \ref{sec-thin-layer}. This fact exemplifies the importance of the boundary conditions in partial differential equations, and it tells us which of the two models should be used when simulating flow of a thin fluid layer between two surfaces: if the fluid pressure is dominant (that is, it is of order $O(\varepsilon^{-2})$), and the fluid velocity is known on the upper and lower surfaces, we must use the lubrication model obtained in section \ref{sec-lubrication};  if the fluid pressure is not dominant (that is, it is of order $O(1)$), and the tractions are known on the upper and lower surfaces, we must use the shallow water model obtained in section \ref{sec-thin-layer}. In the first case we will say that the fluid is ``driven by the pressure'' and in the second that it is ``driven by the velocity''.

In the lubrication model derived in section \ref{sec-lubrication}, the pressure is determined by the equation \eqref{Reynolds_gen}, and it depends on the fluid velocity on the upper and lower surfaces of the domain, and on the speed at which these surfaces move, as well as on the geometry of the surface $\vec{X}$, and on the pressure at $\partial D$ (see remark \ref{rmk-b-c-p}). The fluid velocities inside the domain are subsequently obtained from the pressure using the equations \eqref{u1_0_lub}-\eqref{u2_0_lub}. In the shallow water model of section \ref{sec-thin-layer},  the fluid velocities are calculated from equations \eqref{ec_Vi0-v2} and \eqref{ec_div_eps0ds}, and they are determined by the geometry of the surface $\vec{X}$, as well as by the applied tractions (that is, the pressures $\pi_0^0=\pi_1^0$ and the friction forces), while the fluid pressure is obtained now from the expression \eqref{p0s-v2}.

But, when do we know ``a priori'' if the fluid is ``driven by the pressure'' or ``driven by the velocity'', that is, if we should use the lubrication model or the shallow water model? If we look closely at sections \ref{sec-lubrication} and \ref{sec-thin-layer}, we can say that the lubrication model describes the fluid behavior when the pressure differences at $\partial D$ are large enough, forcing the fluid movement described in \eqref{u1_0_lub}-\eqref{u2_0_lub}, and that the shallow water model describes the fluid behavior when the pressure differences are small at $\partial D$, so that the pressure is determined by the pressure applied to the upper and lower surfaces of the domain and by its separation velocity (see \eqref{p0s-v2}).


\appendix

\appendixpage{}

\section{Change of variable} \label{ApendiceA} 

Let us consider the change of variable \eqref{eq-1-cv}-\eqref{eq-2-cv} between the original domain \eqref{eq-o-domain} and the reference domain \eqref{eq-Omega}. 

Its jacobian matrix is
\begin{displaymath}
    \mathbf{J}^{\varepsilon}=
        \begin{pmatrix}
        \dfrac{ \partial x_1^{\varepsilon}}{\partial \xi_1} &
\dfrac{ \partial x_1^{\varepsilon}}{\partial \xi_2}  & \dfrac{
\partial x_1^{\varepsilon}}{\partial \xi_3}  &
 \dfrac{ \partial x_1^{\varepsilon}}{\partial t}\\
{}\\
        \dfrac{ \partial x_2^{\varepsilon}}{\partial \xi_1}  &
\dfrac{ \partial x_2^{\varepsilon}}{\partial \xi_2}  & \dfrac{
\partial x_2^{\varepsilon}}{\partial \xi_3}  & \dfrac{ \partial
x_2^{\varepsilon}}{\partial t}\\
{}\\
        \dfrac{ \partial x_3^{\varepsilon}}{\partial \xi_1}&
        \dfrac{ \partial x_3^{\varepsilon}}{\partial \xi_2}&
        \dfrac{ \partial x_3^{\varepsilon}}{\partial \xi_3}&
 \dfrac{ \partial x_3^{\varepsilon}}{\partial t}\\
{}\\
        \dfrac{ \partial t^{\varepsilon}}{\partial \xi_1}&
        \dfrac{ \partial t^{\varepsilon}}{\partial \xi_2}&
        \dfrac{ \partial t^{\varepsilon}}{\partial \xi_3}        &
                  \dfrac{ \partial t^{\varepsilon}}{\partial t} \\
        \end{pmatrix}
\end{displaymath}
and it is clear from \eqref{eq-1-cv}-\eqref{base_a3} that 
  \begin{eqnarray}
\dfrac{ \partial x_i^{\varepsilon}}{\partial \xi_j} &=&
a_{ji}+\varepsilon \xi_3 \dfrac{\partial h}{\partial \xi_j} a_{3i}
+\varepsilon \xi_3 h\dfrac{\partial a_{3i}}{\partial \xi_j}, \quad (i=1,2,3; j=1,2) \label{parcial_xiuv_ap}\\
\dfrac{ \partial x_i^{\varepsilon}}{\partial \xi_3} &=& \varepsilon
h
a_{3i}, \quad (i=1,2,3) \label{parcial_xiw_ap}\\
\dfrac{ \partial x_i^{\varepsilon}}{\partial t} &=& \dfrac{\partial
x_i}{\partial t}+\varepsilon \xi_3 \dfrac{\partial h}{\partial
t}a_{3i} +\varepsilon \xi_3 h\dfrac{\partial a_{3i}}{\partial t},
\quad (i=1,2,3) \label{parcial_xit_ap}  \\
\dfrac{ \partial t^{\varepsilon}}{\partial \xi_1} &=&\dfrac{
\partial t^{\varepsilon}}{\partial \xi_2}=\dfrac{ \partial
t^{\varepsilon}}{\partial \xi_3}=0 , 
\label{parcial_tuvw_ap}\\
\dfrac{ \partial t^{\varepsilon}}{\partial t} &=&1
\label{parcial_tt_ap} \end{eqnarray}

We can compute 
\begin{equation}
    (\mathbf{J}^{\varepsilon})^{-1}=
        \begin{pmatrix}
        \dfrac{ \partial \xi_1}{\partial x_1^{\varepsilon}} &
\dfrac{ \partial \xi_1}{\partial x_2^{\varepsilon}}  &
 \dfrac{ \partial \xi_1}{\partial x_3^{\varepsilon}}  &
\dfrac{ \partial \xi_1}{\partial t^{\varepsilon}} \\
{}\\
        \dfrac{ \partial \xi_2}{\partial x_1^{\varepsilon}}  &
\dfrac{ \partial \xi_2}{\partial x_2^{\varepsilon}}  &
 \dfrac{ \partial \xi_2}{\partial x_3^{\varepsilon}}  &
\dfrac{ \partial \xi_2}{\partial t^{\varepsilon}}\\
{}\\
        \dfrac{ \partial \xi_3}{\partial x_1^{\varepsilon}}&
       \dfrac{ \partial \xi_3}{\partial x_2^{\varepsilon}}&
        \dfrac{ \partial \xi_3}{\partial x_3^{\varepsilon}}&
 \dfrac{ \partial \xi_3}{\partial t^{\varepsilon}}\\
{}\\
        \dfrac{ \partial t}{\partial x_1^{\varepsilon}}&
        \dfrac{ \partial t}{\partial x_2^{\varepsilon}}&
        \dfrac{ \partial t}{\partial x_3^{\varepsilon}}         &
        \dfrac{ \partial t}{\partial t^{\varepsilon}} \\
        \end{pmatrix} \label{J-1_ap}
\end{equation}
writing its components in the basis $\left\{
\vec{a}_1,\vec{a}_2,\vec{a}_3\right\}$:
\begin{eqnarray}
\left(\dfrac{ \partial \xi_1}{\partial x_1^{\varepsilon}}, \dfrac{
\partial \xi_1}{\partial x_2^{\varepsilon}}, \dfrac{ \partial \xi_1}{\partial x_3^{\varepsilon}} \right)
&=& \alpha_1 \vec{a}_1 + \beta_1 \vec{a}_2 +\gamma_1 \vec{a}_3
\label{Dxu_base_a_ap}
\\
\left(\dfrac{ \partial \xi_2}{\partial x_1^{\varepsilon}}, \dfrac{
\partial \xi_2}{\partial x_2^{\varepsilon}}, \dfrac{ \partial \xi_2}{\partial x_3^{\varepsilon}} \right)
&=& \alpha_2 \vec{a}_1 + \beta_2 \vec{a}_2 +\gamma_2 \vec{a}_3
\\
\left(\dfrac{ \partial \xi_3}{\partial x_1^{\varepsilon}}, \dfrac{
\partial \xi_3}{\partial x_2^{\varepsilon}}, \dfrac{ \partial \xi_3}{\partial x_3^{\varepsilon}} \right)
&=& \alpha_3 \vec{a}_1 + \beta_3 \vec{a}_2 +\gamma_3 \vec{a}_3
\\
\left(\dfrac{ \partial t}{\partial x_1^{\varepsilon}}, \dfrac{
\partial t}{\partial x_2^{\varepsilon}}, \dfrac{ \partial t}{\partial x_3^{\varepsilon}} \right)
&=& \alpha_4 \vec{a}_1 + \beta_4 \vec{a}_2 +\gamma_4 \vec{a}_3
\label{parcial_uvw_xit_ap}
\end{eqnarray} and using that
\begin{equation}(\mathbf{J}^{\varepsilon})^{-1}\mathbf{J}^{\varepsilon}=I\label{JJ-1I}
\end{equation}

Taking into account that 
\begin{eqnarray} &&\vec{a}_i \cdot
\vec{a}_3=0, \quad (i=1,2),\label{a1a30}\\
&&\|\vec{a}_3\|=1\label{a31}\\
&&\vec{a}_3 \cdot \dfrac{\partial \vec{a}_3}{\partial \xi_i}=0 \quad
(i=1,2)\label{a3da30}
\end{eqnarray}
and introducing the following notation for the  coefficients  of the first and second fundamental forms of the surface parametrized by $\vec{X}$ (here $t$ acts only as a parameter):
\begin{eqnarray} E&=&\vec{a}_1 \cdot\vec{a}_1 \label{coef-E} \\
 F&=&\vec{a}_1 \cdot\vec{a}_2 \label{coef-F} \\
 G&=&\vec{a}_2 \cdot\vec{a}_2 \label{coef-G}
\end{eqnarray}
\begin{eqnarray}
e&=& -\vec{a}_1 \cdot \dfrac{\partial \vec{a}_{3}}{\partial
\xi_1} = \vec{a}_3 \cdot \dfrac{\partial \vec{a}_{1}}{\partial
\xi_1} \label{coef-e} \\
f&=&-
\vec{a}_1 \cdot \dfrac{\partial \vec{a}_{3}}{\partial
\xi_2}=-\vec{a}_2 \cdot\dfrac{\partial \vec{a}_{3}}{\partial
\xi_1}= \vec{a}_3 \cdot \dfrac{\partial \vec{a}_{1}}{\partial
\xi_2}= \vec{a}_3 \cdot \dfrac{\partial \vec{a}_{2}}{\partial
\xi_1} \label{coef-f} \\
g&=&-\vec{a}_2 \cdot\dfrac{\partial \vec{a}_{3}}{\partial
\xi_2}= \vec{a}_3 \cdot \dfrac{\partial \vec{a}_{2}}{\partial
\xi_2} \label{coef-g} 
\end{eqnarray}
we deduce from \eqref{JJ-1I}:
\begin{eqnarray}
\alpha_1 &=&\dfrac{\|\vec{a}_2\|^2+ \varepsilon \xi_3 h \left(\vec{a}_2
\cdot \dfrac{\partial \vec{a}_{3}}{\partial
\xi_2}\right)}{A(\varepsilon)} = \dfrac{G - \varepsilon \xi_3 h g}{A(\varepsilon)}\label{alfa1}\\
\beta_1& =&-\dfrac{ \vec{a}_1 \cdot \vec{a}_{2} +\varepsilon \xi_3 h
\left(\vec{a}_1 \cdot\dfrac{\partial \vec{a}_{3}}{\partial
\xi_2}\right)}{A(\varepsilon)} = -\dfrac{F -\varepsilon \xi_3 h
f}{A(\varepsilon)}\label{beta1}\\
\alpha_2 &=&-\dfrac{\vec{a}_2 \cdot \vec{a}_{1} + \varepsilon \xi_3
h\left( \vec{a}_2 \cdot \dfrac{\partial \vec{a}_{3}}{\partial
\xi_1}\right)}{A(\varepsilon)}=\beta_1 \label{alfa2}\\
\beta_2& =&\dfrac{ \|\vec{a}_1\|^2  +\varepsilon \xi_3 h\left(
\vec{a}_1 \cdot\dfrac{\partial \vec{a}_{3}}{\partial
\xi_1}\right)}{A(\varepsilon)} = \dfrac{ E  -\varepsilon \xi_3 h e}{A(\varepsilon)} \label{beta2}\\
\gamma_i&=&0 \quad (i=1,2) \label{gamma12}\\
\alpha_3 &=&-\dfrac{ \xi_3 }{ h}\left( \alpha_1 \dfrac{\partial h}{\partial \xi_1} + \alpha_2 \dfrac{\partial
h}{\partial \xi_2}
\right)\label{alfa312}\\
\beta_3&=& -\dfrac{ \xi_3 }{ h}\left( \beta_1\dfrac{\partial
h}{\partial \xi_1} + \beta_2 \dfrac{\partial h}{\partial
\xi_2}\right) \label{beta312}\\
\gamma_3&=&\dfrac{1}{\varepsilon h} \label{gamma3}\end{eqnarray}
\begin{eqnarray}
\dfrac{
\partial \xi_1}{\partial t^{\varepsilon}} &=&-(\alpha_1 \vec{a}_1 + \beta_1 \vec{a}_2)\cdot\left(
\dfrac{\partial \vec{X}}{\partial t} + \varepsilon \xi_3 h
\dfrac{\partial \vec{a}_{3}}{\partial t} \right)
\label{parcial_u_teps_alfa_beta_ap}\\
\dfrac{
\partial \xi_2}{\partial t^{\varepsilon}}& =&-( \alpha_2 \vec{a}_1 + \beta_2 \vec{a}_2)\cdot
\left(\dfrac{\partial \vec{X}}{\partial t} + \varepsilon \xi_3 h
\dfrac{\partial \vec{a}_{3}}{\partial t}
\right)\label{parcial_v_teps_alfa_beta_ap}\\
\dfrac{
\partial \xi_3}{\partial t^{\varepsilon}}& =&- (\alpha_3 \vec{a}_1 +\beta_3 \vec{a}_2) \cdot\left(\dfrac{\partial \vec{X}}{\partial t}
+\varepsilon \xi_3 h   \dfrac{\partial \vec{a}_{3}}{\partial t}
\right) -\dfrac{1}{\varepsilon h} \vec{a}_3 \cdot\dfrac{\partial
\vec{X}}{\partial t}- \dfrac{\xi_3}{ h} \dfrac{\partial h}{\partial
t} \label{parcial_w_teps_alfa_beta_ap}\\
\alpha_4 &=&\beta_4=\gamma_4=0\\
\dfrac{\partial t}{\partial x^{\varepsilon}_i}&=&0, \quad (i=1,2,3) \label{dtdxi_0_ap} \\
\dfrac{\partial t}{\partial t^{\varepsilon}}&=&1 \label{dtdteps_1_ap}
\end{eqnarray}
where 
\begin{eqnarray}
A(\varepsilon)&=&\|\vec{a}_1\|^2 \|\vec{a}_2\|^2-\left( \vec{a}_1
\cdot \vec{a}_{2}\right)^2 \nonumber  \\
&+&\varepsilon \xi_3 h\left[\|\vec{a}_2\|^2 \left(\vec{a}_1
\cdot\dfrac{\partial \vec{a}_{3}}{\partial \xi_1}\right) +
\|\vec{a}_1\|^2 \left(\vec{a}_2 \cdot \dfrac{\partial
\vec{a}_{3}}{\partial \xi_2}\right)- \left( \vec{a}_1 \cdot
\vec{a}_{2}\right)\left(\vec{a}_1 \cdot\dfrac{\partial
\vec{a}_{3}}{\partial \xi_2}+\vec{a}_2 \cdot \dfrac{\partial
\vec{a}_{3}}{\partial \xi_1}\right)\right]\nonumber  \\
&+& \varepsilon^2 \xi_3^2 h^2 \left[\left(\vec{a}_1
\cdot\dfrac{\partial \vec{a}_{3}}{\partial \xi_1}
\right)\left(\vec{a}_2 \cdot \dfrac{\partial \vec{a}_{3}}{\partial
\xi_2} \right) -\left( \vec{a}_1 \cdot\dfrac{\partial
\vec{a}_{3}}{\partial \xi_2}\right) \left( \vec{a}_2 \cdot
\dfrac{\partial \vec{a}_{3}}{\partial \xi_1} \right)\right]\nonumber\\
&=&EG-F^2 + \varepsilon \xi_3 h\left(-G e - E g+2 f F \right) +
\varepsilon^2 \xi_3^2 h^2 \left(e g -f^2\right) \label{A}
\end{eqnarray}

If we denote by 
\begin{eqnarray}
A^0&=&\|\vec{a}_1\|^2 \|\vec{a}_2\|^2-\left( \vec{a}_1
\cdot \vec{a}_{2}\right)^2 =EG-F^2=\|\vec{a}_1 \times \vec{a}_2\|^2 \label{A0}  \\
A^1&=&\|\vec{a}_2\|^2 \left(\vec{a}_1 \cdot\dfrac{\partial
\vec{a}_{3}}{\partial \xi_1}\right) + \|\vec{a}_1\|^2
\left(\vec{a}_2 \cdot \dfrac{\partial \vec{a}_{3}}{\partial
\xi_2}\right) \nonumber \\
&&{} - \left( \vec{a}_1 \cdot
\vec{a}_{2}\right)\left(\vec{a}_1 \cdot\dfrac{\partial
\vec{a}_{3}}{\partial \xi_2}+\vec{a}_2 \cdot \dfrac{\partial
\vec{a}_{3}}{\partial \xi_1}\right) = -eG-gE+2fF \label{A^1} \\ A^2&=&\left(\vec{a}_1 \cdot\dfrac{\partial
\vec{a}_{3}}{\partial \xi_1} \right)\left(\vec{a}_2 \cdot
\dfrac{\partial \vec{a}_{3}}{\partial \xi_2} \right) -\left(
\vec{a}_1 \cdot\dfrac{\partial \vec{a}_{3}}{\partial \xi_2}\right)
\left( \vec{a}_2 \cdot \dfrac{\partial \vec{a}_{3}}{\partial \xi_1}
\right)=eg-f^2\label{A2}
\end{eqnarray} 
then we obtain that 
\begin{equation}
A(\varepsilon) = A^0 + \varepsilon \xi_3 h A^1 + \varepsilon^2 \xi_3^2 h^2 A^2 \label{eqAepsilon} 
\end{equation}

We remark (see \cite{Stoker}) that $A^0$, $A^1$ and $A^2$ are related to the Gaussian curvature ($K_G$) of the surface parametrized by $\vec{X}$ and its mean curvature ($K_m$), since 
\begin{eqnarray}
K_G &=& \dfrac{eg-f^2}{EG-F^2} = \dfrac{A^2}{A^0} \label{KG} \\
K_m &=& \dfrac{eG+gE-2fF}{2(EG-F^2)} = - \dfrac{A^1}{2A^0} \label{Km}
\end{eqnarray}

Furthermore, the principal curvatures of $\vec{X}$ are the solutions of the equation 
\begin{equation}
A^0 K_n^2+A^1 K_n+A^2=0
\end{equation}

\section{Coefficients definition} \label{ApendiceB}

In this appendix, we introduce some coefficients that depend only on the lower bound surface parametrization, $\vec{X}$ and other coefficients that depend both on the parametrization and on the gap $h$. We will use these coefficients throughout this article.

In addition to the coefficients that will be defined below, others have been introduced in the body of the paper and in \ref{ApendiceA}: the coefficients of the first and second fundamental forms of the surface parametrized by $\vec{X}$ (denoted by $E, F, G$ and $e,f,g$, respectively), defined in \eqref{coef-E}-\eqref{coef-G} and \eqref{coef-e}-\eqref{coef-g} from the basis $\left\{ \vec{a}_1,\vec{a}_2,\vec{a}_3\right\}$ (see \eqref{base_a1}-\eqref{base_a3}), the coefficients $\alpha_i$, $\beta_i$ and $\gamma_i$ ($i=1,2,3$)  in \eqref{alfa1}-\eqref{gamma3}, and their development in powers of $\varepsilon$ in \eqref{alfaides}-\eqref{beta3n},  $A(\varepsilon)$ and its development in powers of $\varepsilon$ in \eqref{A0}-\eqref{eqAepsilon}, along with its relation with the Gaussian curvature and the mean curvature of the surface parametrized by $\vec{X}$ in \eqref{KG}-\eqref{Km}, and, finally, the definition of $\hat{A}_i^0$ ($i=1,2$) in \eqref{Bi}.

The following coefficients depend only on the parametrization $\vec{X}$: 
\hspace*{-0.5cm}\begin{eqnarray}
B_{11}&=&Ge-Ff \label{B11}\\
B_{12}&=&Gf -Fg\label{B12}\\
B_{21}&=& Ef-Fe \label{B21}\\
B_{22}&=& Eg-Ff \label{B22}\\
C^0_l&=& \alpha_l^0 \left( \vec{a}_1 \cdot \dfrac{\partial \vec{X}}{\partial t}\right) + \beta_l^0 \left( \vec{a}_2 \cdot \dfrac{\partial \vec{X}}{\partial t}\right) \quad (l=1,2) \label{C}\\
H^0_{ilk}&=& \alpha_i^0\left(\vec{a}_1\cdot \dfrac{\partial \vec{a}_k}{\partial \xi_l}\right)+\beta_i^0\left(\vec{a}_2 \cdot \dfrac{\partial \vec{a}_k}{\partial \xi_l}\right)\quad (i,l=1,2; \quad k=1,2,3)\label{H}\\
I&=&\left(\vec{a}_1 \times \dfrac{\partial
\vec{a}_3}{\partial \xi_2}\right)\cdot \vec{a}_3 + \left(\dfrac{\partial
\vec{a}_3}{\partial \xi_1} \times
\vec{a}_2\right)\cdot \vec{a}_3 \label{I}\\
J^0_{lm}&=& \alpha_l^0 \alpha_m^0 E + (\beta_l^0 \alpha_m^0 + \alpha_l^0 \beta_m^0) F + \beta_l^0 \beta_m^0 G \nonumber \\ 
&=&\alpha_l^0 \delta_{m1}+\beta_l^0\delta_{m2} \quad (l,m=1,2) \label{J}\\
{L}^0_{kli} &=&  \sum_{m=1}^2 \left[ \left(  \dfrac{\partial\alpha_l^0}{\partial \xi_m} \delta_{m1}+\dfrac{\partial\beta_l^0}{\partial \xi_m} \delta_{m2}+\alpha_l^0 H^0_{mm1} +\beta_l^0H^0_{mm2}\right) \delta_{ki}\right.\nonumber\\
&+&\left. 2  H^0_{imk} J^0_{lm}\right] \quad (i,l=1,2; \quad k=1,2,3) \label{L}\\
Q^0_{ik}&=&\alpha_i^0 \left( \vec{a}_1 \cdot \dfrac{\partial \vec{a}_k}{\partial t}\right) + \beta_i^0 \left( \vec{a}_2 \cdot \dfrac{\partial \vec{a}_k}{\partial t}\right)-\sum_{l=1}^2 H^0_{ilk}C^0_l  \nonumber \\
&&(i=1,2; \,k=1,2,3) \label{Q} \\
R^0_{ik} &=&  Q^0_{ik}+ H^0_{ik3} \left( \dfrac{\partial \vec{X}}{\partial t} \cdot \vec{a}_3\right) \quad (i=1,2; \,k=1,2) \label{R} 
\end{eqnarray}
\hspace*{-0.5cm}\begin{eqnarray}
S^0_{ik} &=&  \dfrac{I \sqrt{A^0} -A^1}{A^0} \left[ \alpha_i^0\left( \dfrac{\partial \vec{a}_k}{\partial \xi_1} \cdot \vec{a}_3\right) +\beta_i^0 \left( \dfrac{\partial \vec{a}_k}{\partial \xi_2} \cdot \vec{a}_3\right) \right]\nonumber\\
&+& \sum_{m=1}^2 \sum_{l=1}^2 \left[ \left( \alpha^0_i \left( \vec{a}_1 \cdot \dfrac{\partial^2 \vec{a}_k}{\partial \xi_l \partial \xi_m}\right) +  \beta^0_i \left( \vec{a}_2 \cdot \dfrac{\partial^2 \vec{a}_k}{\partial \xi_l \partial \xi_m}\right)\right)J^0_{lm}\right.\nonumber\\
&+&\left.\left(  \dfrac{\partial\alpha_l^0}{\partial \xi_m} \delta_{m1}+\dfrac{\partial\beta_l^0}{\partial \xi_m} \delta_{m2}+\alpha_l^0 H^0_{mm1} +\beta_l^0H^0_{mm2}\right)H^0_{ilk}\right]\nonumber\\
&-&\dfrac{1}{\sqrt{A^0}}\left[ \left(\vec{a}_1 \times \dfrac{\partial
\vec{a}_3}{\partial \xi_2}\right)+ \left(\dfrac{\partial
\vec{a}_3}{\partial \xi_1} \times
\vec{a}_2\right)\right] \cdot \left( \alpha^0_i \dfrac{\partial \vec{a}_k}{\partial \xi_1} + \beta^0_i\dfrac{\partial \vec{a}_k}{\partial \xi_2} \right) \nonumber \\
&&(i,k=1,2)\label{S} 
\end{eqnarray}

\begin{rmk}
Coefficients $B^0_{il}$ and $H^0_{il3}$ are related in the following way:
\begin{equation}
{H^0_{il3}}=-\dfrac{B_{il}}{A^0} \quad (i,l=1,2)
\end{equation}
\end{rmk}

\vspace{-0.2cm}

The following coefficients depend on the parametrization $\vec{X}$ and on function $h$:
\begin{eqnarray}
F^0_i(h)&=&\int_0^1 f^0_i \, d\xi_3+\dfrac{s_0}{\rho h}(\vec{f}^1_{R_1} + \vec{f}^1_{R_0}) \cdot \left( \alpha^0_i \vec{a}_1 + \beta^0_i \vec{a}_2 \right) \quad (i=1,2) \label{F} \\
\vec{\eta}(h)&=&\dfrac{\partial h}{\partial \xi_2} (\vec{a}_1 \times
\vec{a}_3) + h \left(\vec{a}_1 \times \dfrac{\partial
\vec{a}_3}{\partial \xi_2}\right) + \dfrac{\partial h}{\partial
\xi_1} (\vec{a}_3\times \vec{a}_2) + h \left(\dfrac{\partial
\vec{a}_3}{\partial \xi_1} \times
\vec{a}_2\right)\label{eta}\\
\psi(h)_{ijl}^0&=& \dfrac{1}{h} \left[  
\left( \alpha_l^0 \dfrac{\partial h}{\partial \xi_1} +\beta_l^0 \dfrac{\partial h}{\partial \xi_2} \right) \delta_{ij} + \dfrac{\partial h}{\partial \xi_j} \left(  \alpha_l^0 \delta_{1i} + \beta_l^0 \delta_{2i}\right)\right] \quad (i,j,l=1,2) \label{psi}\\
\chi(h)_{ik}^0&=&\dfrac{1}{h}   \left\{  \dfrac{\partial h}{\partial
\xi_1} \left[\sum_{l=1}^2H^0_{ilk}\alpha_l^0 -\dfrac{1}{\sqrt{A^0}}(\vec{a}_3\times \vec{a}_2)\cdot \left( \alpha_i^0
\dfrac{
\partial \vec{a}_{k}}{\partial \xi_1}+\beta_i^0 \dfrac{
\partial \vec{a}_{k}}{\partial \xi_2}\right)\right]\right.\nonumber\\
&+&\left. \dfrac{\partial h}{\partial \xi_2} \left[ \sum_{l=1}^ 2H^0_{ilk}\beta_l^0 -\dfrac{1}{\sqrt{A^0}} (\vec{a}_1 \times
\vec{a}_3) \cdot \left( \alpha_i^0
\dfrac{
\partial \vec{a}_{k}}{\partial \xi_1}+ \beta_i^0
\dfrac{
\partial \vec{a}_{k}}{\partial \xi_2} \right)  \right]\right\} \nonumber\\
&& (i=1,2, \,k=1,2,3)\label{chiik}\\
\kappa(h)^0_i&=& -\dfrac{1}{h^2}\left[\alpha^0_i \dfrac{\partial}{\partial t}\left(h \dfrac{
\partial h}{\partial \xi_1}  
\right) +\beta^0_i \dfrac{\partial}{\partial t}\left(h \dfrac{
\partial h}{\partial \xi_2}  
\right) \right]\nonumber\\
&+& \left[ \dfrac{\partial }{\partial \xi_1}\left( \dfrac{\partial \vec{X}}{\partial t} \cdot \vec{a}_3\right) \left( L^0_{31i} - \dfrac{A^1}{A^0} \alpha^0_i \right)+\dfrac{\partial }{\partial \xi_2}\left( \dfrac{\partial \vec{X}}{\partial t} \cdot \vec{a}_3\right) \left( L^0_{32i} - \dfrac{A^1}{A^0} \beta^0_i \right)\right]\nonumber\\
&+& \left( \dfrac{\partial \vec{X}}{\partial t} \cdot \vec{a}_3\right) \left\{ \chi(h)^0_{i3} + \sum_{m=1}^2 \sum_{l=1}^2 \left[ \left( \alpha^0_i \left( \vec{a}_1 \cdot \dfrac{\partial^2 \vec{a}_3}{\partial \xi_l \partial \xi_m}\right) +  \beta^0_i \left( \vec{a}_2 \cdot \dfrac{\partial^2 \vec{a}_3}{\partial \xi_l \partial \xi_m}\right)\right)J^0_{lm}\right.\right.\nonumber\\
&+&\left. \left.\left(  \dfrac{\partial\alpha_l^0}{\partial \xi_m} \delta_{m1}+\dfrac{\partial\beta_l^0}{\partial \xi_m} \delta_{m2}+\alpha_l^0 H^0_{mm1} +\beta_l^0H^0_{mm2}\right)H^0_{il3}\right]\right.\nonumber\\
&-&\left.\dfrac{1}{\sqrt{A^0}}\left[ \left(\vec{a}_1 \times \dfrac{\partial
\vec{a}_3}{\partial \xi_2}\right)+ \left(\dfrac{\partial
\vec{a}_3}{\partial \xi_1} \times
\vec{a}_2\right)\right] \cdot \left( \alpha^0_i \dfrac{\partial \vec{a}_3}{\partial \xi_1} + \beta^0_i\dfrac{\partial \vec{a}_3}{\partial \xi_2} \right)\right\}  \quad (i=1,2)\label{kappa} \\
\hat{\kappa}(h)^0_i &=& \kappa(h)^0_i - \left( \alpha_i^0 \dfrac{
\partial}{\partial \xi_1} \left ( \frac{2}{h} \dfrac{\partial h}{\partial t} \right ) + \beta_i^0 \dfrac{
\partial}{\partial \xi_2} \left ( \frac{2}{h} \dfrac{\partial h}{\partial t} \right ) \right) \quad (i=1,2)\label{kappa-hat}
\end{eqnarray}
where $\delta_{ij}$ is the Kronecker Delta.

\section{Derivation of equations to calculate $\vec{V}^0$} \label{ApendiceC}

Let us identify the terms of order $\varepsilon^0$ in the equation \eqref{ec_ns_ij_alfa_beta}. We simplify that equation, taking into account \eqref{parcial_xi1_teps_alfa_beta}-\eqref{parcial_xi3_teps_alfa_beta}, \eqref{u_3^0}, \eqref{u_i^0s} and \eqref{u_1^1s}-\eqref{p0s-v2}. Then we multiply the equation obtained by $a_{1i}$ and we yield:
\begin{eqnarray}
&&\hspace*{-0.5cm}
\sum_{k=1}^3 \left(\dfrac{ \partial V_k^0}{\partial t} (\vec{a}_1 \cdot \vec{a}_{k}) +
V_k^0 \left(\vec{a}_1 \cdot \dfrac{
\partial \vec{a}_{k}}{\partial t} \right)\right) \nonumber\\
&&{}+ \sum_{l=1}^2 \sum_{k=1}^3\left( \dfrac{
\partial V_k^0}{\partial \xi_l} (\vec{a}_1 \cdot \vec{a}_{k})+ V_k^0   \left(\vec{a}_1 \cdot \dfrac{
\partial \vec{a}_{k}}{\partial \xi_l} \right) \right)\left( V_l^0-C^0_l \right) \nonumber\\
&&=-\dfrac{1}{\rho_0}  \dfrac{
\partial p^0}{\partial \xi_1}  + \nu \left\{ \sum_{m=1}^2 \sum_{l=1}^2 \sum_{k=1}^3
\dfrac{\partial}{\partial \xi_m}\left[ \dfrac{
\partial (V_k^0 \vec{a}_k)}{\partial \xi_l} (\alpha_l^0
{a}_{1j} + \beta_l^0 {a}_{2j}) \right]\cdot \vec{a}_1  (\alpha_m^0
{a}_{1j} + \beta_m^0 {a}_{2j})\right.\nonumber\\
&& \left.{} +\dfrac{1 }{h} \dfrac{A^1}{A^0}  \sum_{k=1}^2 \dfrac{
\partial u_k^1}{\partial \xi_3} (\vec{a}_1 \cdot \vec{a}_{k})
+\dfrac{1 }{h^2} \sum_{k=1}^2  \dfrac{
\partial^2 u_k^2}{\partial \xi_3^3}  (\vec{a}_1 \cdot \vec{a}_{k}) \right\}+ \sum_{k=1}^2
f_k^0 (\vec{a}_1 \cdot \vec{a}_{k}) \label{ec_ns_eps0_abrev_a1}\end{eqnarray}
where we have denoted by $V_3^0=u_3^0$ (to achieve a more compact expression), and coefficients $C^0_l$, $(l=1,2)$, are given by \eqref{C}.

Analogously, if we multiply the same equation by $a_{2i}$ and $a_{3i}$, we obtain, respectively:
\begin{eqnarray}
&&\hspace*{-0.5cm}
\sum_{k=1}^3 \left(\dfrac{ \partial V_k^0}{\partial t} (\vec{a}_2 \cdot \vec{a}_{k}) +
V_k^0 \left(\vec{a}_2 \cdot \dfrac{
\partial \vec{a}_{k}}{\partial t} \right)\right) \nonumber\\
&&{}+ \sum_{l=1}^2 \sum_{k=1}^3 \left( \dfrac{
\partial V_k^0}{\partial \xi_l} (\vec{a}_2 \cdot \vec{a}_{k})+ V_k^0   \left(\vec{a}_2 \cdot \dfrac{
\partial \vec{a}_{k}}{\partial \xi_l} \right) \right)\left( V_l^0- C^0_l \right) \nonumber\\
&&=-\dfrac{1}{\rho_0}  \dfrac{
\partial p^0}{\partial \xi_2}  + \nu \left\{ \sum_{m=1}^2 \sum_{l=1}^2 \sum_{k=1}^3
\dfrac{\partial}{\partial \xi_m}\left[ \dfrac{
\partial (V_k^0 \vec{a}_k)}{\partial \xi_l} (\alpha_l^0
{a}_{1j} + \beta_l^0 {a}_{2j}) \right]\cdot \vec{a}_2  (\alpha_m^0
{a}_{1j} + \beta_m^0 {a}_{2j})\right.\nonumber\\
&& \left.{} +\dfrac{1 }{h} \dfrac{A^1}{A^0}  \sum_{k=1}^2 \dfrac{
\partial u_k^1}{\partial \xi_3} (\vec{a}_2 \cdot \vec{a}_{k})
+\dfrac{1 }{h^2} \sum_{k=1}^2  \dfrac{
\partial^2 u_k^2}{\partial \xi_3^3}  (\vec{a}_2 \cdot \vec{a}_{k}) \right\}+ \sum_{k=1}^2
f_k^0 (\vec{a}_2 \cdot \vec{a}_{k}) \label{ec_ns_eps0_abrev_a2}\end{eqnarray}
\begin{eqnarray}
&&\hspace*{-0.5cm}\dfrac{ \partial V_3^0}{\partial t}+
\sum_{k=1}^2 
V_k^0 \left(\vec{a}_3 \cdot \dfrac{
\partial \vec{a}_{k}}{\partial t} \right)\nonumber\\
&&{}+ \sum_{l=1}^2  \left( \dfrac{
\partial V_3^0}{\partial \xi_l}+ \sum_{k=1}^2 V_k^0   \left(\vec{a}_3 \cdot \dfrac{
\partial \vec{a}_{k}}{\partial \xi_l} \right) \right)\left( V_l^0- C^0_l\right) \nonumber\\
&&=-\dfrac{1}{\rho_0 h}  \dfrac{
\partial p^1}{\partial \xi_3}  + \nu \left\{ \sum_{m=1}^2 \sum_{l=1}^2 \sum_{k=1}^3
\dfrac{\partial}{\partial \xi_m}\left[ \dfrac{
\partial (V_k^0 \vec{a}_k)}{\partial \xi_l} (\alpha_l^0
{a}_{1j} + \beta_l^0 {a}_{2j}) \right]\cdot \vec{a}_3  (\alpha_m^0
{a}_{1j} + \beta_m^0 {a}_{2j})\right.\nonumber\\
&& \left.{} +\dfrac{1 }{h} \dfrac{A^1}{A^0}  \dfrac{\partial h}{\partial t}
+\dfrac{1 }{h^2}   \dfrac{
\partial^2 u_3^2}{\partial \xi_3^3}  \right\}+ 
f_3^0 \label{ec_ns_eps0_abrev_a3}\end{eqnarray}

Next we multiply equation \eqref{ec_ns_eps0_abrev_a1} by $\alpha_1^0$ and we add equation \eqref{ec_ns_eps0_abrev_a2} multiplied by $\alpha_2^0$ to get:
\begin{eqnarray}
&&\hspace*{-0.5cm}\dfrac{ \partial V_1^0}{\partial t}   +
\sum_{l=1}^2 \alpha_l^0
\sum_{k=1}^3 \left(
V_k^0 \left(\vec{a}_l \cdot \dfrac{
\partial \vec{a}_{k}}{\partial t} \right)\right) \nonumber\\
&&{}+ \sum_{l=1}^2 \left( \dfrac{
\partial V_1^0}{\partial \xi_l} + \sum_{m=1}^2 \alpha_m^0 \sum_{k=1}^3  V_k^0   \left(\vec{a}_m \cdot \dfrac{
\partial \vec{a}_{k}}{\partial \xi_l} \right) \right)\left( V_l^0- C^0_l \right)=-\dfrac{1}{\rho_0}  \sum_{l=1}^2 \alpha_l^0 \dfrac{
\partial p^0}{\partial \xi_l} \nonumber\\
&&{} + \nu \left\{ \sum_{m=1}^2\sum_{p=1}^2\alpha_p^0 \sum_{l=1}^2 \sum_{k=1}^3
\dfrac{\partial}{\partial \xi_m}\left[ \dfrac{
\partial (V_k^0 \vec{a}_k)}{\partial \xi_l} (\alpha_l^0
{a}_{1j} + \beta_l^0 {a}_{2j}) \right]\cdot \vec{a}_p  (\alpha_m^0
{a}_{1j} + \beta_m^0 {a}_{2j})\right.\nonumber\\
&& \left.{} +\dfrac{1 }{h} \dfrac{A^1}{A^0}   \dfrac{
\partial u_1^1}{\partial \xi_3}
+\dfrac{1 }{h^2}  \dfrac{
\partial^2 u_1^2}{\partial \xi_3^3}   \right\}+ 
f_1^0  \label{ec_ns_eps0_abrev_a1b}\end{eqnarray}

In the same way, we multiply equation \eqref{ec_ns_eps0_abrev_a2} by $\beta_2^0$ and we add equation \eqref{ec_ns_eps0_abrev_a1} multiplied by $\beta_1^0$ to obtain:
\begin{eqnarray}
&&\hspace*{-0.5cm}\dfrac{ \partial V_2^0}{\partial t}   +
\sum_{l=1}^2 \beta_l^0
\sum_{k=1}^3 \left(
V_k^0 \left(\vec{a}_l \cdot \dfrac{
\partial \vec{a}_{k}}{\partial t} \right)\right) \nonumber\\
&&{}+ \sum_{l=1}^2 \left( \dfrac{
\partial V_2^0}{\partial \xi_l} + \sum_{m=1}^2 \beta_m^0 \sum_{k=1}^3  V_k^0   \left(\vec{a}_m \cdot \dfrac{
\partial \vec{a}_{k}}{\partial \xi_l} \right) \right)\left( V_l^0-C^0_l\right) =-\dfrac{1}{\rho_0}  \sum_{l=1}^2 \beta_l^0 \dfrac{
\partial p^0}{\partial \xi_l} \nonumber\\
&&{} + \nu \left\{ \sum_{m=1}^2\sum_{p=1}^2\beta_p^0 \sum_{l=1}^2 \sum_{k=1}^3
\dfrac{\partial}{\partial \xi_m}\left[ \dfrac{
\partial (V_k^0 \vec{a}_k)}{\partial \xi_l} (\alpha_l^0
{a}_{1j} + \beta_l^0 {a}_{2j}) \right]\cdot \vec{a}_p  (\alpha_m^0
{a}_{1j} + \beta_m^0 {a}_{2j})\right.\nonumber\\
&& \left.{} +\dfrac{1 }{h} \dfrac{A^1}{A^0}   \dfrac{
\partial u_2^1}{\partial \xi_3}
+\dfrac{1 }{h^2}  \dfrac{
\partial^2 u_2^2}{\partial \xi_3^3}   \right\}+ 
f_2^0  \label{ec_ns_eps0_abrev_a2b}\end{eqnarray}

We yield the following equations by integrating  \eqref{ec_ns_eps0_abrev_a1b}-\eqref{ec_ns_eps0_abrev_a2b} 
over  $\xi_3$ from 0 to 1, and using expressions \eqref{u_1^1s}-\eqref{u_2^1s},  \eqref{V_1^1_V^0_abrev}-\eqref{V_2^1_V^0_abrev} and 
\eqref{dif_du_1^2_dxi3_1_0}-\eqref{dif_du_2^2_dxi3_1_0}:

\begin{eqnarray}
&&\hspace*{-0.5cm}\dfrac{ \partial V_1^0}{\partial t}  + \sum_{l=1}^2 \left( V_l^0-C^0_l \right) \dfrac{
\partial V_1^0}{\partial \xi_l} +
\sum_{k=1}^3 
\left( Q^0_{1k}+\sum_{l=1}^2
H^0_{1lk} V_l^0\right) V_k^0 \nonumber\\
&&=-\dfrac{1}{\rho_0}  \sum_{l=1}^2 \alpha_l^0 \dfrac{
\partial p^0}{\partial \xi_l}+\nu \left\{ \sum_{m=1}^2 \sum_{l=1}^2 
 \dfrac{
\partial^2 V_1^0 }{\partial \xi_m \partial \xi_l} J^0_{lm}
+ \sum_{k=1}^2 \sum_{l=1}^2  \dfrac{
\partial V_k^0 }{\partial \xi_l}( L^0_{kl1} +\psi(h)^0_{1kl} )\right.\nonumber\\
&& \left. {}
+ \sum_{k=1}^2 V_k^0 ({S}_{1k}^0+\chi(h)^0_{1k})  + \kappa(h)^0_1  \right\}+ {F}^0_1(h) \label{ec_V10s} \end{eqnarray}
\begin{eqnarray}
&&\hspace*{-0.5cm}\dfrac{ \partial V_2^0}{\partial t}+ \sum_{l=1}^2 \left( V_l^0-C^0_l \right)  \dfrac{
\partial V_2^0}{\partial \xi_l}  +\sum_{k=1}^3 \left(Q^0_{2k} + \sum_{l=1}^2  H^0_{2lk} V_l^0 \right)
V_k^0\nonumber\\
&&=-\dfrac{1}{\rho_0} \sum_{l=1}^2 \beta_l^0 \dfrac{
\partial p^0}{\partial \xi_l}   + \nu \left\{  \sum_{m=1}^2 \sum_{l=1}^2  \dfrac{
\partial^2 V_2^0}{\partial \xi_m \partial \xi_l} J^0_{lm} +
\sum_{k=1}^2 \sum_{l=1}^2  \dfrac{
\partial V_k^0 }{\partial \xi_l}( L^0_{kl2} +\psi(h)^0_{2kl} )
 \right.\nonumber\\
&& \left.{} +\sum_{k=1}^2 V_k^0 ({S}^0_{2k}+\chi(h)^0_{2k}) + \kappa(h)^0_2
\right\} + {F}^0_2(h)  \label{ec_V20s}\end{eqnarray}
where $H^0_{ilk}$, $J^0_{lm}$, $L^0_{kli}$, $Q^0_{ik}$, $S^0_{ik}$, $F^0_i(h)$, $\psi(h)^0_{ikl}$, $\chi(h)^0_{ik}$ and $\kappa(h)^0_i$ are given by \eqref{H}, \eqref{J}-\eqref{Q}, 
\eqref{S}-\eqref{F} and \eqref{psi}-\eqref{kappa}.

Finally, from last equations, taking into account that $V_3^0 = u_3^0$, \eqref{u_3^0}, \eqref{p0s-v2} and rearranging terms, we obtain 
\begin{eqnarray}
&&\hspace*{-0.5cm}\dfrac{ \partial V_i^0}{\partial t}  + \sum_{l=1}^2 \left( V_l^0-C^0_l \right) \dfrac{
\partial V_i^0}{\partial \xi_l} +
\sum_{k=1}^2 
\left( R^0_{ik}+\sum_{l=1}^2
H^0_{ilk} V_l^0 \right) V_k^0\nonumber\\
&&=-\dfrac{1}{\rho_0}\left( \alpha_i^0 \dfrac{
\partial \pi_0^0}{\partial \xi_1}  + \beta_i^0 \dfrac{
\partial \pi_0^0}{\partial \xi_2} \right) \nonumber \\ 
&&{} +\nu \left\{ \sum_{m=1}^2 \sum_{l=1}^2 
 \dfrac{
\partial^2 V_i^0 }{\partial \xi_m \partial \xi_l} J^0_{lm}
+ \sum_{k=1}^2 \sum_{l=1}^2  \dfrac{
\partial V_k^0 }{\partial \xi_l}( L^0_{kli} +\psi(h)^0_{ikl} )\right.\nonumber\\
&& \left. {}
+ \sum_{k=1}^2 V_k^0 ({S}_{ik}^0+\chi(h)^0_{ik})  + \hat{\kappa}(h)^0_i  \right\}+ {F}^0_i(h)-
 Q^0_{i3} \left ( \frac{\partial \vec{X}}{\partial t}
\cdot \vec{a}_3\right )  \quad (i=1,2)\label{ec_Vi0} \end{eqnarray}
where $R^0_{ik}$ and $\hat{\kappa}(h)^0_i$ are given by \eqref{R} and \eqref{kappa-hat}.

\end{document}